\documentclass{siamltex}

\usepackage{amsmath, amssymb}

\usepackage{hyperref}
\usepackage{bm}

\usepackage{enumitem}

\usepackage{graphicx}
\usepackage{epstopdf}

\usepackage{makecell}

\usepackage{xcolor}
\usepackage{colortbl}

\usepackage{hhline}
\usepackage{multirow}
\usepackage{tabu}

\usepackage[linesnumbered,ruled]{algorithm2e}
\usepackage{setspace}

\colorlet{lightgray}{gray!12}
\colorlet{darkgray}{black!68}

\DeclareMathOperator*{\argavg}{arg\ avg}

\DeclareMathOperator*{\essinf}{ess\ inf}

\usepackage{epstopdf}

\usepackage{subfigure}

\usepackage{tikz}

\usetikzlibrary{shapes.misc}
\tikzset{cross/.style={cross out, draw=black, minimum size=2*(#1-\pgflinewidth), inner sep=0pt, outer sep=0pt},
cross/.default={4pt}}

\newtheorem{thm}{Theorem}[section]
\newtheorem{examp}{Example}[section]

\newtheorem{rem}{Remark}[section]

\newtheorem{assump}{Assumption}[section]

\usepackage{amsfonts}
\usepackage{amssymb,amsmath,bm}
\usepackage{graphicx,color,rotating,fancybox}

\usepackage{hyperref}

\begin{document}
\bibliographystyle{plain}

\pagestyle{myheadings}
\markboth{A.J.~Crowder, C.E.~Powell, A.~Bespalov}{}

\title{{Efficient adaptive multilevel stochastic Galerkin approximation using implicit a posteriori error estimation}\thanks{{This work was supported 
by  EPSRC grants EP/P013317/1 and EP/P013791/1. The second author would like to thank the Isaac Newton Institute for Mathematical Sciences, Cambridge, for support and hospitality during the Uncertainty Quantification programme as well as the Simons Foundation. This work was partially also supported by EPSRC grant no EP/K032208/1.}}}

\author{A.J.~Crowder\thanks{School of Mathematics, University of Manchester, 
Oxford Road, Manchester M13 9PL, United Kingdom (adam.crowder@postgrad.manchester.ac.uk).}
\and
C.E.~Powell\thanks{School of Mathematics, 
University of Manchester, Oxford Road, Manchester M13 9PL, 
United Kingdom (c.powell@manchester.ac.uk).}
\and
A.~Bespalov\thanks{School of Mathematics, 
University of Birmingham, Edgbaston, Birmingham B15 2TT, 
United Kingdom (a.bespalov@bham.ac.uk).}
}
\maketitle

\begin{abstract} Partial differential equations (PDEs) with inputs that depend on infinitely many parameters pose serious theoretical and computational challenges.  Sophisticated numerical algorithms that automatically determine which parameters need to be activated in the approximation space in order to estimate a quantity of interest to a prescribed error tolerance are needed. For elliptic PDEs with parameter-dependent coefficients, stochastic Galerkin finite element methods (SGFEMs) have been well studied. Under certain assumptions, it can be shown that there exists a sequence of SGFEM approximation spaces for which the energy norm of the error decays to zero at a rate that is independent of the number of input parameters. However, it is not clear how to adaptively construct these spaces in a practical and computationally efficient way.  We present a new adaptive SGFEM algorithm that tackles elliptic PDEs with parameter-dependent coefficients quickly and efficiently. We consider approximation spaces with a multilevel structure---where each solution mode is associated with a finite element space on a potentially different mesh---and use an implicit a posteriori error estimation strategy to steer the adaptive enrichment of the space. At each step, the components of the error estimator are used to assess the potential benefits of a variety of enrichment strategies, including whether or not to activate more parameters. No marking or tuning parameters are required. Numerical experiments for a selection of test problems demonstrate that the new method performs optimally in that it generates a sequence of approximations for which the estimated energy error decays to zero at the same rate as the error for the underlying finite element method applied to the associated parameter-free problem.
\end{abstract}



\section{Introduction}\label{Sec:intro}

In many engineering and other real world applications, we frequently encounter models consisting of partial differential equations (PDEs) which have uncertain or parameter-dependent inputs. When the solutions are sufficiently smooth with respect to these parameters, it is known that stochastic Galerkin finite element methods (SGFEMs) \cite{MR1083354,MR1870425,MR2084236}, also known as intrusive polynomial chaos methods in the statistics and engineering communities, offer a powerful alternative to brute force sampling methods for propagating uncertainty to the model outputs. When the number of input parameters in the PDE model is countably infinite (which may arise, for example, if we represent an uncertain spatially varying coefficient as a Karhunen-Lo\`eve expansion), then we encounter significant theoretical and numerical challenges. In general, it is not known a priori which parameters need to be incorporated into discretisations of the model in order to estimate specific quantities of interest to a prescribed error tolerance. Ad hoc selection of a finite subset of parameters prior to applying a standard SGFEM is computationally convenient, but may lead to inaccurate results with no guaranteed error bounds. In this work we consider the steady-state diffusion problem with a spatially varying coefficient that depends on infinitely many parameters, and develop a computationally efficient multilevel SGFEM which uses an a posteriori error estimator to adaptively construct appropriate approximation spaces. 

Let the spatial domain $D \subset \mathbb{R}^{2}$ be bounded with a Lipschitz polygonal boundary $\partial D$ and let $y_{1}, y_{2}, \ldots$ be a countable sequence of parameters with $y_{m} \in \Gamma_{m} = [-1,1],$ for $m \in \mathbb{N}$. We consider the parametric diffusion problem: find $u(\mathbf{x},\mathbf{y}):D\times\Gamma\to\mathbb{R}$ that satisfies
\begin{align}
-\nabla\cdot(a(\mathbf{x},\mathbf{y})\nabla u(\mathbf{x},\mathbf{y})) &= f(\mathbf{x}),& &\mathbf{x}\in D,\ \mathbf{y}\in \Gamma,\label{intro_PDE_problem1}\\
u(\mathbf{x},\mathbf{y})&=0,& & \mathbf{x}\in\partial D,\ \mathbf{y} \in\Gamma. \label{intro_PDE_problem2}
\end{align}
Here, $\mathbf{y} = [y_1,y_2,\dots]^{\top}\in\Gamma$ where $\Gamma= \Pi_{m=1}^{\infty}\Gamma_{m}$ is the parameter domain. The coefficient $a(\mathbf{x}, \mathbf{y})$ should be positive and bounded on $D \times \Gamma$.  We also make the following important assumption.
\smallskip
\begin{assump}\label{Assump:a_decomp}\normalfont The coefficient $a(\mathbf{x},\mathbf{y})$ admits the decomposition
\begin{align}\label{a_decomp}
a(\mathbf{x},\mathbf{y}) = a_0(\mathbf{x}) + \sum_{m=1}^\infty a_m(\mathbf{x})y_m,
\end{align}
with $a_0(\mathbf{x}),\ a_m(\mathbf{x})\in L^\infty(D)$ and $ || a_m||_{L^\infty(D)} \to 0$ sufficiently quickly as $m\to\infty$ so that
\begin{align}\label{a_conv}
\sum_{m=1}^\infty||a_m||_{L^\infty(D)}< \essinf_{\mathbf{x}\in D} a_0(\mathbf{x}).
\end{align}
\end{assump}\noindent
Note that \eqref{a_conv} helps to ensure the well-posedness of the weak formulation of \eqref{intro_PDE_problem1}--\eqref{intro_PDE_problem2}. This will be made more rigorous in the next section. 

Standard SGFEMs seek approximations to $u(\mathbf{x},\mathbf{y})$ in \eqref{intro_PDE_problem1}--\eqref{intro_PDE_problem2} in a tensor product space $X$ of the form
\begin{align}\label{single_level_X}
X:= H_1\otimes P,\qquad H_1 := \text{span}\{\phi_i(\mathbf{x})\}_{i=1}^n,\qquad P:=\text{span}\{\psi_j(\mathbf{y})\}_{j=1}^s,
\end{align}
where $H_1$ is a finite element space associated with a mesh ${\cal T}_{h}$ on the spatial domain $D$ and $P$ is a set of polynomials on the parameter domain $\Gamma$ in a \emph{finite} number (say, $M$) of the parameters $y_{m}$. In this case, $u_X\in X$ admits the decomposition 
$$u_X(\mathbf{x}, \mathbf{y}) = \sum_{j=1}^s u_j(\mathbf{x}) \psi_j(\mathbf{y}), \qquad u_{j} \in H_{1}.$$ 
We use the term `single-level'  approximation to mean that $X$ is defined as in \eqref{single_level_X}. Here, each coefficient $u_j$ is associated with the \emph{same} finite element space $H_{1}$. In contrast, we will work with spaces $X$ which have a `multilevel' structure, by which we mean that the coefficients $u_j$ may each reside in a \emph{different} finite element space. These finite element spaces will be associated with a sequence of meshes which each have a different `level' number.


Handling inputs of the form \eqref{a_decomp} is a non-trivial task. Suppose we truncate $a(\mathbf{x},\mathbf{y})$ in \eqref{a_decomp} after $M$ terms (assuming that $||a_{m} ||_{\infty}  \ge ||a_{m+1} ||_{\infty}$) and define $X$ as in \eqref{single_level_X}, where $\mathbf{y}=[y_{1}, \ldots, y_{M}]^{\top}$.  A priori error estimates provided in \cite{MR2084236} reveal that the rate of convergence of standard SGFEMs deteriorates as $M \to \infty$. This phenomenon is referred to as the \emph{curse of dimensionality}. Many recent works provide a priori error analysis for more sophisticated SGFEMs in the case where we have infinitely many parameters. For example, see \cite{MR2317004,MR2500242,MR2566594,MR2728424,MR2763359,MR3091365,MR3298364}. In each of these works, the decay rate, or equivalently, the summability of the sequence $\{\|a_m\|_{\infty}\}_{m=1}^\infty$ plays an important role.  Various theoretical results have been established proving the existence of a sequence of SGFEM approximation spaces $X^{0}, X^{1}, \ldots$, such that the energy norm of the error decays to zero at a rate that is independent of the number of parameters, as $N_{\textrm{dof}} = \textrm{dim}(X) \to \infty$. These results all assume that $X$ has a more complex structure than in \eqref{single_level_X} but demonstrate that SGFEMs \emph{can} be immune to the curse of dimensionality if implemented in the right way. 

In \cite{MR2728424,MR2763359,MR3091365} a multilevel structure is imposed on $X$. Theoretical results show that if $ \|a_m\|_{\infty} \to 0 $ fast enough, then there exists a sequence of multilevel spaces for which the error decays to zero at the rate afforded to the chosen finite element method for the parameter-free analogue of \eqref{intro_PDE_problem1}--\eqref{intro_PDE_problem2}. Given a sequence of finite element spaces (with different level numbers), we use an implicit a posteriori error estimation scheme to design an appriopriate sequence of multilevel SGFEM spaces. By implicit, we mean that the approach uses the residual associated with the SGFEM solution indirectly and requires the solution of additional problems. Starting with an initial low-dimensional space $X^{0}$, the resulting energy error is estimated. The components of the error estimator are then examined to steer the enrichment of $X^{0}$.  Adaptive schemes have also been proposed in \cite{MR3154028,MR3042573,MR3423228,MR3578034}, but using an explicit error estimation strategy which uses the residual directly.  Explicit error estimators often lead to less favourable effectivity indices than implicit schemes.  Moreover, the algorithms presented in \cite{MR3154028,MR3042573,MR3423228,MR3578034} all rely on a D\"{o}rfler-like marking strategy \cite{MR1393904}, and require the selection of multiple tuning or marking parameters. The optimal selection of these is unclear, however, and is problem-dependent. The authors of \cite{MR3177362,MR3519560,MR3396480,MR3771406} consider single-level approximation spaces and implement an implicit error estimation strategy. We revisit \cite{MR3177362,MR3519560}, extend the error estimation strategy considered there to the more complex multilevel setting, and use this to design an accurate and efficient adaptive multilevel SGFEM algorithm.

\subsection{Outline} In Section \ref{sec2} we introduce the weak formulation of \eqref{intro_PDE_problem1}--\eqref{intro_PDE_problem2} and review conditions for well-posedness. In Section \ref{Sec:SGFEM_approx} we describe the multilevel construction of SGFEM approximation spaces and give practical information about how to assemble the matrices associated with the discrete problem in a computationally efficient way.  In Section \ref{Sec:apost} we extend the implicit energy norm a posteriori error estimation strategy developed in \cite{MR3177362,MR3519560} for SGFEM approximation spaces $X$ of the form \eqref{single_level_X} to the multilevel setting. In Section \ref{Sec:5} we introduce a new adaptive  algorithm that uses the error estimation strategy from Section \ref{Sec:apost} to design problem-dependent multilevel SGFEM approximation spaces. Numerical results are presented in Section \ref{Sec:num_results}.


\section{Weak Formulation of the Parametric Diffusion Problem}\label{sec2} We assume that $y_m\in \Gamma_m := [-1,1]$ for each $m\in \mathbb{N}$ and that $\pi_m$ is a measure on $(\Gamma_m,\mathcal{B}(\Gamma_m))$, where $\mathcal{B}(\Gamma_m)$ denotes the Borel $\sigma$--algebra on $\Gamma_m$. We also assume that
\begin{align}\label{mean_zero_rv}
\int_{\Gamma_m}y_m\ d\pi_m(y_m) = 0,\quad m\in\mathbb{N}.
\end{align}
For instance, this is true when $y_m$ is the image of a mean zero random variable and $\pi_m$ is the associated probability measure. We assume that $y_{m}$ is the image of a uniform random variable $\xi_{m} \sim U([-1,1])$ and so the associated probability measure $\pi_{m}$ has density $\rho_{m}=1/2$ with respect to Lebesgue measure.  We now define the parameter domain $\Gamma = \Pi_{m=1}^{\infty} \Gamma_{m}$ and the product measure 
$$\pi(\mathbf{y}):= \prod_{m=1}^\infty \pi_m(y_m).$$
If the parameters $y_{m}$ are images of independent random variables then the associated probability measure has this separable form.

We are interested in Galerkin approximations of $u$ satisfying \eqref{intro_PDE_problem1}--\eqref{intro_PDE_problem2} and thus start by considering its variational formulation:
\begin{align}\label{para_weak_prob}
\text{find } u\in  V:= L_\pi^2(\Gamma, H_0^1(D)) :\quad B(u,v) = F(v),\quad\text{for all }v \in V.
\end{align}
Here, $H_0^1(D)$ is the usual Hilbert space of functions that vanish on $\partial D$ in the sense of trace and $L_\pi^2(\Gamma)$ is the space of functions that are  square integrable with respect to $\pi(\mathbf{y})$ on $\Gamma$. That is,
\begin{align*}
L_\pi^2(\Gamma) := \bigg\{v(\mathbf{y})\ |\ \langle v,v\rangle_{L^{2}_{\pi}(\Gamma)} = \int_\Gamma v(\mathbf{y})^2\ d\pi(\mathbf{y}) < \infty \bigg\}.
\end{align*}
The space $V$ is equipped with the norm $||\cdot||_V$, where
\begin{align*}
||v||_V = \left(\int_\Gamma ||v(\cdot,\mathbf{y})||_{H_0^1(D)}^2\ d\pi(\mathbf{y})\right)^{\frac{1}{2}},
\end{align*}
and $||v||_{H_0^1(D)} = ||\nabla v||_{L^2(D)}$ for all $v\in H_0^1(D)$. The bilinear form $B:V\times V\to\mathbb{R}$ and the linear functional $F:V\to\mathbb{R}$ are defined by
\begin{align}
B(u,v) &= \int_\Gamma\int_D a(\mathbf{x},\mathbf{y})\nabla u(\mathbf{x},\mathbf{y})\cdot \nabla v(\mathbf{x},\mathbf{y})\ d\mathbf{x}\ d\pi(\mathbf{y}),\label{B_bilinear_form}\\
F(v) &= \int_\Gamma\int_D f(\mathbf{x})v(\mathbf{x},\mathbf{y})\ d\mathbf{x}\ d\pi(\mathbf{y}).
\end{align}
To ensure that \eqref{para_weak_prob} is well-posed, $B(\cdot,\cdot)$ must be bounded and coercive over $V$. This is ensured by the following assumption.
\begin{assump}\label{Assump:a_bounds}\normalfont There exist real positive constants $a_{\min}$ and $a_{\max}$ such that
\begin{align*}
0 < a_{\min} \leq a(\mathbf{x},\mathbf{y}) \leq a_{\max} < \infty,\quad a.e. \text{ in } D\times\Gamma.
\end{align*}
\end{assump}\noindent
Note that \eqref{a_conv} is a sufficient condition for Assumption \ref{Assump:a_bounds} to hold. If Assumption \ref{Assump:a_bounds} holds, the bilinear form \eqref{B_bilinear_form} induces a norm (the so-called energy norm),
 $$||v||_B = B(v,v)^{1/2}, \qquad \textrm{ for all } v \in V.$$
In addition, to ensure that $F(\cdot)$ is bounded on $V$ we assume $f(\mathbf{x})\in L^2(D)$. We will also make the following assumption.
\begin{assump}\label{Assump:a0_bounds}\normalfont There exist real positive constants $a_{\min}^0$ and $a_{\max}^0$ such that
\begin{align*}
0 < a_{\min}^0 \leq a_0(\mathbf{x}) \leq a_{\max}^0 < \infty,\quad a.e. \text{ in } D.
\end{align*}
\end{assump}

Due to \eqref{a_decomp},  we have the decomposition,
\begin{align}\label{B_decomposition}
B(u,v) = B_0(u,v) + \sum_{m=1}^\infty B_m(u,v),\quad\text{for all }u,v\in V,
\end{align}
where the component bilinear forms are given by
\begin{align}
B_0(u,v) &= \int_\Gamma\int_D a_0(\mathbf{x})\nabla u(\mathbf{x},\mathbf{y})\cdot \nabla v(\mathbf{x},\mathbf{y})\ d\mathbf{x}\ d\pi(\mathbf{y}),\label{B0_bilinear_form}\\
B_m(u,v) &= \int_\Gamma\int_D a_my_m(\mathbf{x})\nabla u(\mathbf{x},\mathbf{y})\cdot \nabla v(\mathbf{x},\mathbf{y})\ d\mathbf{x}\ d\pi(\mathbf{y}).\label{Bm_bilinear_form}
\end{align}
If Assumption \ref{Assump:a0_bounds} holds, the bilinear form \eqref{B0_bilinear_form} also induces the norm $||v||_{B_0} = B_0(v,v)^{1/2}$ on $V$, associated with the coefficient $a_{0}$. It is then straightforward to show that
\begin{align*}
\lambda ||v||_B^2 \leq ||v||_{B_0}^2 \leq \Lambda||v||_B^2,\quad\text{for all } v\in V,
\end{align*}
where $0< \lambda < 1 < \Lambda < \infty$ and
\begin{align}\label{lambdas_def}
\lambda:=a_{\min}^0a_{\max}^{-1},\qquad \lambda:=a_{\max}^0a_{\min}^{-1},
\end{align}
and so the norms $|| \cdot ||_{B}$ and $|| \cdot ||_{B_{0}}$ are equivalent.

\section{Multilevel SGFEM Approximation}\label{Sec:SGFEM_approx}
We can compute a Galerkin approximation to $u\in V$ by projecting \eqref{para_weak_prob} onto a finite-dimensional subspace $X\subset V$. The best known rates of convergence with respect to $N_{\textrm{dof}}=\textrm{dim}(X)$  (see \cite{MR3298364, MR2728424,MR2763359,MR3091365}) are achieved for approximation spaces that have a multilevel structure, which we now describe. As usual, we exploit the fact that $V\cong H_0^1(D)\otimes L_\pi^2(\Gamma)$ and construct $X$ by tensorising separate subspaces of  $H_0^1(D)$ and $L_\pi^2(\Gamma)$.

For the parameter domain, we first introduce families of univariate polynomials $\{\psi_{n}(y_m)\}_{n\in\mathbb{N}_0}$ on $\Gamma_m$ for each $m=1,2,\dots$ that are orthonormal with respect to the inner product 
$$\langle v, w \rangle_{L_{\pi_m}^2(\Gamma_m)} = \int_{\Gamma_{m}} v(y_{m}) w(y_{m}) d \pi_{m}(y_{m}).$$
 Here, $n$ denotes the polynomial degree and $\psi_0(y_m) = 1$.  Now we define the set of finitely supported multi-indices $J:= \{ \mu = (\mu_1,\mu_2,\dots)\in \mathbb{N}_0^{\mathbb{N}};\ \#\text{supp}(\mu) < \infty \}$ where $\text{supp}(\mu):= \{m\in\mathbb{N};\ \mu_m\neq 0\}$ and consider multivariate tensor product polynomials of the form
\begin{align}\label{psi_def_construction}
\psi_\mu(\mathbf{y}) = \prod_{m=1}^\infty \psi_{\mu_m}(y_m) = \prod_{m\in\text{supp}(\mu)}\psi_{\mu_m}(y_m),\quad \mu\in J.
\end{align}
The countable set $\{\psi_\mu(\mathbf{y})\}_{\mu\in J}$ is an orthonormal basis of $L_\pi^2(\Gamma)$ with respect to the inner product $\langle \cdot,\cdot \rangle_{L_\pi^2(\Gamma)}$. Orthonormality comes from the separability of $\pi(\mathbf{y})$ and the construction \eqref{psi_def_construction} since
\begin{align}\label{mutual_orthogonality}
\langle\psi_\mu(\mathbf{y}),\psi_\nu(\mathbf{y})\rangle_{L_\pi^2(\Gamma)} = \prod_{m=1}^\infty \langle \psi_{\mu_m}(y_m),\psi_{\nu_m}(y_m)\rangle_{L_{\pi_m}^2(\Gamma_m)} = \prod_{m=1}^\infty \delta_{\mu_m\nu_m} = \delta_{\mu\nu},
\end{align}
for all $\mu,\nu\in J$. Now, given any finite set $J_P\subset J$ (which we assume always contains the multi-index $\mu =  (0,0,\dots)$) we can construct a finite-dimensional set $P := \left\{\psi_\mu(\mathbf{y}),   \mu\in J_P \right\}\subset L_\pi^2(\Gamma)$ of multivariate polynomials on $\Gamma$. Note that we can also write 
\begin{align*}
P = \bigoplus_{\mu\in J_P} P^{\mu},\quad P^{\mu} =  \text{span}\{\psi_\mu(\mathbf{y})\},\quad \mu\in J_P.
\end{align*}

Given a set of multi-indices $J_P$, we will construct approximation spaces of the form
\begin{align}\label{multilevel_X}
X := \bigoplus_{\mu\in J_P} X^\mu:= \bigoplus_{\mu\in J_P} H_1^{\mu} \otimes P^{\mu} \subset V,
\end{align}
where each $H_1^{\mu}\subset H_0^1(D)$ is a finite element space associated with the spatial domain $D$ and
\begin{align*}
H_1^{\mu} := \text{span}\bigg\{\phi_i^{\mu}(\mathbf{x});\ i=1,2,\dots,N_1^{\mu}\bigg\},\quad\text{for all } \mu\in J_P.
\end{align*}
For each $\mu\in J_P$ we may use a potentially different space $H_1^\mu$. Compare $X$ in \eqref{multilevel_X} to $X$ in \eqref{single_level_X}. The latter can be written as $X := \oplus_{\mu\in J_P} H_1 \otimes P^{\mu}$.  To work with spaces of the form \eqref{multilevel_X}, we need to select an appropriate  set  $\mathbf{H}_1 := \{H_1^\mu\}_{\mu\in J_P}$ of finite element spaces. To this end, we assume that we can construct a nested sequence of meshes ${\cal T}_{i}$, $i=0,1, \ldots$ (of rectangular or triangular elements) that give rise to a sequence of conforming finite element spaces $H^{(0)} \subset H^{(1)}\subset \cdots H^{(i)}\cdots \subset H_{0}^{1}(D)$. In this setting, the index $i$ denotes the mesh `level number'. We will assume that the polynomial degree is fixed in the definition of the finite element spaces, and only the mesh is changing as we change the level. If $j>i$, then ${\cal T}_{j}$ can be obtained from ${\cal T}_{i}$ by one or more mesh refinements.

For notational convenience, we collect the meshes into a set
\begin{align}\label{sequence_of_meshes}
\bm{\mathcal{T}} := \big\{\mathcal{T}_{i};\ i=0,1,2,\dots\big\}.
\end{align}
For each $\mu \in J_{P}$, the space $H_1^\mu$ is constructed using one of the meshes from $\bm{\mathcal{T}}$. That is, to each $\mu \in J_{P}$ we assign a mesh level number $\ell^\mu=i$ (for some $i\in\mathbb{N}_0$) and set $H_1^\mu = H^{(i)}$. If $\ell^\mu = \ell^\nu$ for some $\mu,\nu\in J_P$, then $H_1^\mu=H_1^\nu$. We collect the chosen levels $\ell^{\mu}$ in the set $\bm\ell := \{\ell^{\mu}\}_{\mu\in J_P}$. Now, any space $X$ of the form \eqref{multilevel_X} is determined by choosing a finite set $J_P$ of multi-indices and a set $\bm\ell$ of associated mesh level numbers. Clearly, $\text{card}(\bm\ell) = \text{card}(J_P) < \infty$. 

Once $J_{P}$ and $\bm\ell$ have been chosen, our SGFEM approximation $u_X\in X$ to $u \in V$ is found by solving the discrete problem:
\begin{align}\label{para_discrete_prob}
\text{find } u_X\in X:\quad B(u_X,v) = F(v),\quad\text{for all }v \in X.
\end{align}
For $u_X$ to be computable, it is essential that the sum in \eqref{B_decomposition} has a finite number of nonzero terms. Let $M\in \mathbb{N}$ be the smallest integer such that $\mu_{m} = 0$ for all $m > M$ and for all $\mu\in J_P$. That is, let $M$ be the number of parameters $y_m$ that are `active' in the definition of $J_P$. Then, provided \eqref{mean_zero_rv} holds, $B_m(u_X,v) = 0$ for $u_X,v\in X$ for all $m>M$ (e.g. see \cite{MR3177362}). In other words, the choice of $J_P$ implicitly truncates the sum after $M$ terms; we do not have to truncate $a(\mathbf{x}, \mathbf{y})$ a priori. Expanding the Galerkin approximation as
\begin{align}\label{uX_expansion}
u_X = \sum_{\mu\in J_P}u_X^{\mu}(\mathbf{x})\psi_\mu(\mathbf{y}),\qquad u_X^{\mu} = \sum_{i=1}^{N_1^{\mu}} u_i^{\mu}\phi_i^{\mu}(\mathbf{x}),\qquad u_i^{\mu}\in\mathbb{R},\end{align}
and taking test functions $v = \psi_\nu(\mathbf{y})\phi_j^{\nu}(\mathbf{x})$ for all $\nu\in J_P$ and $j = 1,2,\dots,N_1^{\nu}$ yields a system of $N_{\textrm{dof}}$ equations $A\mathbf{u} = \mathbf{b}$ for the unknown coefficients $u_i^\mu$ that define $u_X$, where
$$N_{\textrm{dof}} = \sum_{\mu \in J_{P}} \textrm{dim}(X^{\mu}) = \sum_{\mu \in J_{P}} N_{1}^{\mu}.$$
If multilevel SGFEMs are to be useful in practice, we have to be able to assemble the components of this linear system and solve it efficiently. We discuss this next.

\subsection{Multilevel SGFEM Matrices}
The matrix $A$ and the vectors $\mathbf{b}$ and $\mathbf{u}$ each have a block structure, with the blocks indexed by the elements (multi-indices) of $J_P$, namely
\begin{align*}
[A_{\mu\nu}]_{ij} = [A_{\nu\mu}]_{ji} &= B\big(\psi_\mu\phi_i^{\mu},\psi_\nu\phi_j^{\nu}\big)\qquad (A\textrm{ is symmetric}),\\
[\mathbf{b}_\nu]_j &= F\big(\psi_\nu\phi_j^{\nu}\big),\\
[\mathbf{u}_\mu]_i &= u_i^{\mu},
\end{align*}
for $i= 1,2,\dots,N_1^{\mu}$ and $j = 1,2,\dots,N_1^{\nu}$. For single-level methods, the resulting system matrix admits the Kronecker product structure (e.g., see \cite{MR2491431}) $ K_0\otimes G_0 + \sum_{m=1}^M K_m\otimes G_m,$ where $\left\{K_m\right\}_{m=0}^{M}$ are stiffness matrices associated with the \emph{same} finite element space and 
\begin{align*}
[G_0]_{\mu\nu} = [G_0]_{\nu\mu} = \delta_{\nu\mu},\qquad  [G_m]_{\mu\nu}=[G_m]_{\nu\mu} = \int_\Gamma y_m\psi_\mu(\mathbf{y})\psi_\nu(\mathbf{y})\ d\pi(\mathbf{y}),\qquad m=1,2,\dots,M.
\end{align*}
In the multilevel approach, there is no such Kronecker structure. The $\nu\mu$\textsuperscript{th} block of $A$ is given by
\begin{align}\label{mult_A_struc}
A_{\nu\mu} = \sum_{m=0}^M [G_m]_{\nu\mu} K_{\nu\mu}^m,\qquad [K_{\nu\mu}^m]_{ji} = \int_D a_m(\mathbf{x})\nabla\phi_i^{\mu}(\mathbf{x})\cdot \nabla \phi_j^{\nu}(\mathbf{x})\ d\mathbf{x},
\end{align}
for $i=1,2,\dots,N_1^{\mu}$ and $j = 1,2,\dots,N_1^{\nu}$. The entries of the stiffness matrix $K_{\nu\mu}^m$ in \eqref{mult_A_struc} depend on basis functions associated with a pair of meshes $\mathcal{T}_{\ell^{\mu}}$ and $\mathcal{T}_{\ell^{\nu}}$, which may be different. Consequently, $K_{\nu\mu}^m$ is non-square if $\ell^{\mu}\neq \ell^{\nu}$ for any $\mu,\nu\in J_P$.  


The key to a fast and efficient multilevel SGFEM algorithm is to first determine what, and what does not, need computing.  If we use iterative solvers, then we only need to compute the action of $A$ on vectors. Here, $\mathbf{v} = A\mathbf{x}$ can be computed blockwise via
\begin{align}\label{multi_mat_vec}
[\mathbf{v}]_\nu = [A\mathbf{x}]_\nu = \sum_{\mu\in J_P}A_{\nu\mu} [\mathbf{x}]_\mu = \sum_{\mu\in J_P} \sum_{m=0}^M [G_m]_{\nu\mu} K_{\nu\mu}^m [\mathbf{x}]_\mu,\qquad \nu\in J_P.
\end{align}
We need only compute $K_{\nu\mu}^m $ for all \emph{distinct} triplets $(m,\ell^{\nu},\ell^{\mu})$ where the corresponding entry $[G_m]_{\nu\mu}$ is \emph{non-zero}.  Due to the orthonormality of the polynomials $\{\psi_\mu(\mathbf{y})\}_{\mu\in J_P}$, the matrices $\{G_m\}_{m=0}^M$ are very sparse (in fact $G_{0}=I$). Indeed, if the density $\rho_{m}$ associated with $\pi_m$ on $\Gamma_m$ is an even function (symmetric about zero), then the matrices $\{G_m\}_{m=1}^M$ have at most two nonzero entries per row, see \cite{MR2491431,MR2644740}. Hence, a naive upper bound for the number of required stiffness matrices is $(1 + 2M)\textrm{card}(J_P)$. This takes the sparsity of $G_{m}$ into account, but does not exploit the fact that the \emph{same} mesh may be assigned to several multi-indices $\mu \in J_{P}$. An adaptive algorithm for automatically selecting $J_{P}$ and the associated set of mesh level numbers $\bm\ell$ is developed in Section \ref{Sec:5}. In Table \ref{Tab:Kcomp} we record $\textrm{card}(J_{P})$ and the number of matrices $K_{\nu\mu}^m$ that are required at the final step of that algorithm (when the error tolerance is set to $\epsilon=2 \times 10^{-3}$), for the test problems outlined in Section \ref{Sec:num_results} (see also Table \ref{Tab:exp1_1}). Since the same mesh level number is assigned to many multi-indices in $J_{P}$, the number of matrices computed is significantly lower than the bound.

\begin{table}[t!]
\centering
\caption{Naive upper bound for the number of matrices $K_{\nu\mu}^m$ that need computing for the test problems (TP.1--TP.4) outlined in Section \ref{Sec:num_results}, and the actual number required. The set $J_{P}$ and the mesh level numbers $\bm\ell$ are selected automatically using Algorithm 1 in Section \ref{Sec:5}. See Sections \ref{Sec:exp_setup} and \ref{Sec:exp_1} for more details.}\label{Tab:Kcomp}
\setlength{\tabcolsep}{3pt}
\begin{tabu}{|[1.3pt] c |[1.3pt]c|c|c|c|[1.3pt]}
\tabucline[1.3pt]{-}
Test Problem & $\textrm{card}(J_P)$ & $M$  & $(1 + 2M)\textrm{card}(J_P)$ & actual \\
\hline
TP.1 & 169 &  93 & 31,603 & 616\\
TP.2 & 36  &  13 &    972 &  96\\
TP.3 & 17 &   3  &    119 &  35\\
TP.4 & 21 &   8  &    357 &  54\\
\tabucline[1.3pt]{-}
\end{tabu}
\end{table}

Adaptive multilevel SGFEMs have been considered in \cite{MR3042573,MR3154028}. Those works use an explicit a posteriori error estimation strategy to drive the enrichment of the approximation space. In \cite{MR3154028}, all stiffness matrices $K_{\nu\mu}^m$ that are non-square ($\ell^\nu\neq \ell^\mu$) are approximated using a projection technique involving only the square matrices $K_{\mu\mu}^m$ that feature in the diagonal blocks $A_{\mu\mu}$ of $A$. Even with this approximation, the multilevel approach considered in \cite{MR3154028} is reported to be computationally expensive. In the next section, we describe how the matrices $K_{\nu\mu}^m$ can be computed quickly and efficiently, without the need for the approximation used in \cite{MR3154028}. 

\subsection{Assembly of Stiffness Matrices}
We describe the construction of $K_{\nu\mu}^m$ for two multi-indices $\mu,\nu\in J_P$, with $\ell^{\mu}\neq \ell^{\nu}$ ($m$ is not important here) for a simple example. For clarity of presentation, we consider uniform meshes of square elements. However, the procedure is applicable to any conforming FEM spaces $H_1^\mu$ and $H_1^\nu$ for which $\mathcal{T}_{\ell^\nu}$ is nested in $\mathcal{T}_{\ell^\mu}$, or equivalently, when $\mathcal{T}_{\ell^\nu}$ is obtained from a conforming (without introducing hanging nodes) refinement of $\mathcal{T}_{\ell^\mu}$.

\smallskip
\begin{examp}\label{Examp:stiff_construct}\normalfont
For simplicity, assume that $D\subset\mathbb{R}^2$ is a square and $H_1^{\mu}$ and $H_1^{\nu}$ are spaces of continuous piecewise bilinear functions associated with two uniform meshes of square elements ($\mathbb{Q}_1$ elements). In particular, let $\mathcal{T}_{\ell^\mu}$ denote a uniform $2\times 2$ square partition of $D$ with mesh level number $\ell^\mu$ and let $\mathcal{T}_{\ell^\nu}$ be a uniform $4\times 4$ square partition of $D$ with $\ell^\nu:=\ell^\mu+1$ (representing, in this case, a uniform refinement of $\mathcal{T}_{\ell^\mu}$). For now, we retain the boundary nodes so that $N_1^{\mu} := \textrm{dim}(H_1^\mu) = 9$ and $N_1^{\nu} := \textrm{dim}(H_1^\nu) = 25$. See Figures \ref{Fig:incompatible_meshes_a} and \ref{Fig:incompatible_meshes_b}. To construct $K_{\nu\mu}^m\in\mathbb{R}^{25 \times 9}$, we compute a \emph{coarse-element} matrix for each element  $\Box_{\text{coarse}}$ in $\mathcal{T}_{\ell^{\mu}}$. In Figure \ref{Fig:incompatible_meshes_c} we highlight one such element, and the four (fine) elements ${\Box_{\text{fine}}}$ in $\mathcal{T}_{\ell^{\nu}}$ that are embedded within it. The associated coarse-element matrix $K_{\nu\mu,c}^m\in\mathbb{R}^{9\times 4}$ has entries
\begin{align*}
[K_{\nu\mu,c}^m]_{ji} = \int_{\Box_{\text{coarse}}} a_m(\mathbf{x})\nabla\phi_i^{\mu,c}(\mathbf{x})\cdot \nabla \phi_j^{\nu,c}(\mathbf{x})\ d\mathbf{x},\qquad i = 1,2,3,4,\qquad j = 1,2,\dots,9,
\end{align*}
where $\{\phi_i^{\mu,c}\}_{i=1}^4$ and $\{\phi_j^{\nu,c}\}_{j=1}^9$  are basis functions associated with the round and cross markers, with support on $\Box_{\text{coarse}}$ and patches of $\Box_{\text{coarse}}$, respectively.
To construct $K_{\nu\mu,c}^m$, we concatenate four \emph{fine-element} matrices $K_{\nu\mu,f}^m\in\mathbb{R}^{4\times4}$ defined by
\begin{align*}
[K_{\nu\mu,f}^m]_{ji} = \int_{\Box_{\text{fine}}} a_m(\mathbf{x})\nabla\phi_i^{\mu,c}(\mathbf{x})\cdot \nabla \phi_j^{\nu,f}(\mathbf{x})\ d\mathbf{x},\qquad i,j = 1,2,3,4,
\end{align*}
where $\Box_{\text{fine}}$ is one of the four elements embedded in $\Box_{\text{coarse}}$. Here, $\{\phi_j^{\nu,f}\}_{j=1}^4$ are the basis functions defined with respect to the crosses in Figure \ref{Fig:local_matrix_patches}, that are supported only on $\Box_{\text{fine}}$ (shaded region). 

For $\mathbb{Q}_{1}$ elements,  constructing $K_{\nu\mu}^m$ boils down to the assembly of $4\times 4$ fine-element matrices $K_{\nu\mu,f}^m$. Similarly, for $\mathbb{Q}_2$ elements (continuous piecewise biquadratic approximation), the procedure requires the assembly of $9\times 9$ fine-element matrices $K_{\nu\mu,f}^m$. If $K_{\nu\mu}^m$ is square ($\ell^\mu = \ell^\nu$), we can use the traditional element construction. In either case, we only need to perform integration on elements in the fine mesh, for which we implement an exact quadrature rule.
\end{examp}

\begin{figure}[t!]
\centering
\subfigure[$\mathcal{T}_{\ell^{\mu}}$.]{
\begin{minipage}[c]{0.25\linewidth}
\centering
\begin{tikzpicture}
\fill[lightgray] (0,0) rectangle (2.5,2.5);
\draw[step=1.25cm,darkgray,line width = 0.6mm] (0,0) grid (2.5,2.5);
\filldraw[fill=white, draw = black, very thick](0   ,0   ) circle (0.12cm);
\filldraw[fill=white, draw = black, very thick](1.25,0   ) circle (0.12cm);
\filldraw[fill=white, draw = black, very thick](2.5 ,0   ) circle (0.12cm);
\filldraw[fill=white, draw = black, very thick](0   ,1.25) circle (0.12cm);
\filldraw[fill=white, draw = black, very thick](1.25,1.25) circle (0.12cm);
\filldraw[fill=white, draw = black, very thick](2.5 ,1.25) circle (0.12cm);
\filldraw[fill=white, draw = black, very thick](2.5 ,0   ) circle (0.12cm);
\filldraw[fill=white, draw = black, very thick](0   ,2.5 ) circle (0.12cm);
\filldraw[fill=white, draw = black, very thick](1.25,2.5 ) circle (0.12cm);
\filldraw[fill=white, draw = black, very thick](2.5 ,2.5 ) circle (0.12cm);
\end{tikzpicture}
\vspace{5pt}
\end{minipage}
\label{Fig:incompatible_meshes_a}
}
\subfigure[$\mathcal{T}_{\ell^{\nu}}$.]{
\begin{minipage}[c]{0.25\linewidth}
\centering
\begin{tikzpicture}
\fill[lightgray] (0,0) rectangle (2.5,2.5);
\draw[step=0.625,darkgray,thick] (0,0) grid (2.5,2.5);
\draw[very thick] (0    ,0    ) node[cross] {};
\draw[very thick] (0.625,0    ) node[cross] {};
\draw[very thick] (1.25 ,0    ) node[cross] {};
\draw[very thick] (1.875,0    ) node[cross] {};
\draw[very thick] (2.5  ,0    ) node[cross] {};
\draw[very thick] (0    ,0.625) node[cross] {};
\draw[very thick] (0.625,0.625) node[cross] {};
\draw[very thick] (1.25 ,0.625) node[cross] {};
\draw[very thick] (1.875,0.625) node[cross] {};
\draw[very thick] (2.5  ,0.625) node[cross] {};
\draw[very thick] (0    ,1.25 ) node[cross] {};
\draw[very thick] (0.625,1.25 ) node[cross] {};
\draw[very thick] (1.25 ,1.25 ) node[cross] {};
\draw[very thick] (1.875,1.25 ) node[cross] {};
\draw[very thick] (2.5  ,1.25 ) node[cross] {};
\draw[very thick] (0    ,1.875) node[cross] {};
\draw[very thick] (0.625,1.875) node[cross] {};
\draw[very thick] (1.25 ,1.875) node[cross] {};
\draw[very thick] (1.875,1.875) node[cross] {};
\draw[very thick] (2.5  ,1.875) node[cross] {};
\draw[very thick] (0    ,2.5  ) node[cross] {};
\draw[very thick] (0.625,2.5  ) node[cross] {};
\draw[very thick] (1.25 ,2.5  ) node[cross] {};
\draw[very thick] (1.875,2.5  ) node[cross] {};
\draw[very thick] (2.5  ,2.5  ) node[cross] {};
\end{tikzpicture}
\end{minipage}
\label{Fig:incompatible_meshes_b}
}
\subfigure[$\Box_{\text{coarse}}$ and its embedded elements.]{
\begin{minipage}[c]{0.25\linewidth}
\centering
\begin{tikzpicture}
\fill[lightgray] (0,0) rectangle (1.25,1.25);
\draw[step=0.625,thick,darkgray] (0,0) grid (2.5,2.5);
\draw[step=1.25cm,line width = 0.6mm,darkgray] (0,0) grid (2.5,2.5);
\filldraw[fill=white, draw = black, very thick](0  ,0  ) circle (0.2cm);
\filldraw[fill=white, draw = black, very thick](1.25,0  ) circle (0.2cm);
\filldraw[fill=white, draw = black, very thick](0  ,1.25) circle (0.2cm);
\filldraw[fill=white, draw = black, very thick](1.25,1.25) circle (0.2cm);
\draw[very thick] (0    ,0    ) node[cross] {};
\draw[very thick] (0.625,0    ) node[cross] {};
\draw[very thick] (1.25 ,0    ) node[cross] {};
\draw[very thick] (0    ,0.625) node[cross] {};
\draw[very thick] (0.625,0.625) node[cross] {};
\draw[very thick] (1.25 ,0.625) node[cross] {};
\draw[very thick] (0    ,1.25 ) node[cross] {};
\draw[very thick] (0.625,1.25 ) node[cross] {};
\draw[very thick] (1.25 ,1.25 ) node[cross] {};
\end{tikzpicture}
\end{minipage}
\label{Fig:incompatible_meshes_c}
}
\caption{Example meshes with (a) $N_1^{\mu} = 9$ and level number $\ell^{\mu}$ and (b) $N_1^{\nu} = 25$ and level number $\ell^{\nu} = \ell^{\mu}+1$.
}
\label{Fig:incompatible_meshes}
\end{figure}
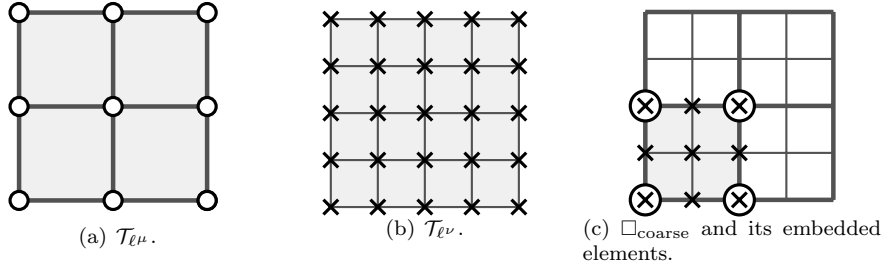

\begin{figure}[t!]
\centering
\subfigure{
\begin{minipage}[c]{0.2\linewidth}
\centering
\begin{tikzpicture}
\fill[lightgray] (0,0) rectangle (1,1);
\draw[step=1cm,darkgray,thick,dashed] (0,0) grid (2,2);
\draw[step=1cm,darkgray,line width = 0.6mm] (0,0) grid (1,1);
\filldraw[fill=white, draw = black, very thick](0   ,0   ) circle (0.2cm);
\filldraw[fill=white, draw = black, very thick](2   ,0   ) circle (0.12cm);
\filldraw[fill=white, draw = black, very thick](0   ,2   ) circle (0.12cm);
\filldraw[fill=white, draw = black, very thick](2   ,2   ) circle (0.12cm);
\draw[very thick] (0  ,0  ) node[cross] {};
\draw[very thick] (1  ,0  ) node[cross] {};
\draw[very thick] (0  ,1  ) node[cross] {};
\draw[very thick] (1  ,1  ) node[cross] {};
\end{tikzpicture}
\end{minipage}
}
\subfigure{
\begin{minipage}[c]{0.2\linewidth}
\centering
\begin{tikzpicture}
\fill[lightgray] (1,0) rectangle (2,1);
\draw[step=1cm,darkgray,thick,dashed] (0,0) grid (2,2);
\draw[step=1cm,darkgray,line width = 0.6mm] (1,0) grid (2,1);
\filldraw[fill=white, draw = black, very thick](0   ,0   ) circle (0.12cm);
\filldraw[fill=white, draw = black, very thick](2   ,0   ) circle (0.2cm);
\filldraw[fill=white, draw = black, very thick](0   ,2   ) circle (0.12cm);
\filldraw[fill=white, draw = black, very thick](2   ,2   ) circle (0.12cm);
\draw[very thick] (2  ,0  ) node[cross] {};
\draw[very thick] (1  ,0  ) node[cross] {};
\draw[very thick] (2  ,1  ) node[cross] {};
\draw[very thick] (1  ,1  ) node[cross] {};
\end{tikzpicture}
\end{minipage}
}
\subfigure{
\begin{minipage}[c]{0.2\linewidth}
\centering
\begin{tikzpicture}
\fill[lightgray] (0,1) rectangle (1,2);
\draw[step=1cm,darkgray,thick,dashed] (0,0) grid (2,2);
\draw[step=1cm,darkgray,line width = 0.6mm] (0,1) grid (1,2);
\filldraw[fill=white, draw = black, very thick](0   ,0   ) circle (0.12cm);
\filldraw[fill=white, draw = black, very thick](2   ,0   ) circle (0.12cm);
\filldraw[fill=white, draw = black, very thick](0   ,2   ) circle (0.2cm);
\filldraw[fill=white, draw = black, very thick](2   ,2   ) circle (0.12cm);
\draw[very thick] (0  ,2  ) node[cross] {};
\draw[very thick] (1  ,2  ) node[cross] {};
\draw[very thick] (0  ,1  ) node[cross] {};
\draw[very thick] (1  ,1  ) node[cross] {};
\end{tikzpicture}
\end{minipage}
}
\subfigure{
\begin{minipage}[c]{0.2\linewidth}
\centering
\begin{tikzpicture}
\fill[lightgray] (1,1) rectangle (2,2);
\draw[step=1cm,darkgray,thick,dashed] (0,0) grid (2,2);
\draw[step=1cm,darkgray,line width = 0.6mm] (1,1) grid (2,2);
\filldraw[fill=white, draw = black, very thick](0   ,0   ) circle (0.12cm);
\filldraw[fill=white, draw = black, very thick](2   ,0   ) circle (0.12cm);
\filldraw[fill=white, draw = black, very thick](0   ,2   ) circle (0.12cm);
\filldraw[fill=white, draw = black, very thick](2   ,2   ) circle (0.2cm);
\draw[very thick] (2  ,2  ) node[cross] {};
\draw[very thick] (1  ,2  ) node[cross] {};
\draw[very thick] (2  ,1  ) node[cross] {};
\draw[very thick] (1  ,1  ) node[cross] {};
\end{tikzpicture}
\end{minipage}
}
\caption{The four embedded elements in Figure \ref{Fig:incompatible_meshes_c} on which we construct four $4\times 4$ local matrices.}
\label{Fig:local_matrix_patches}
\end{figure}
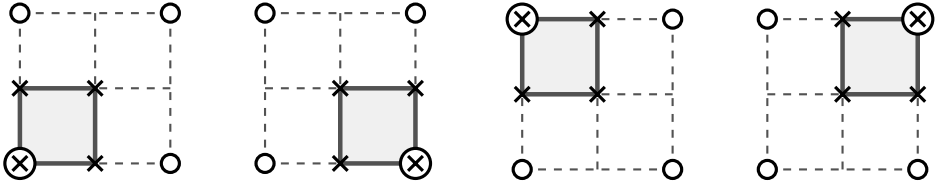

\smallskip
\begin{rem} When the meshes $\mathcal{T}_{\ell^\mu}$ and $\mathcal{T}_{\ell^\nu}$ are uniform, as in Example \ref{Examp:stiff_construct}, the computation of the fine-element matrices can be vectorised over all the coarse elements.
\end{rem}

\section{Energy Norm A Posteriori Error Estimation}\label{Sec:apost}
Given an approximation space $X$ of the form \eqref{multilevel_X} and an SGFEM approximation $u_X\in X$ satisfying \eqref{para_discrete_prob}, we want to estimate the energy error $||u-u_X||_B$. We now extend the implicit strategy developed in \cite{MR3177362,MR3519560}.

Computing the error $e=u-u_X\in V$ is a non-trivial task. Due to the bilinearity of $B(\cdot,\cdot)$ it is clear that $e$ satisfies
\begin{align*}
B(e,v) = B(u,v) - B(u_X,v) = F(v) - B(u_X,v),\quad\text{for all } v\in V.
\end{align*}
We look for an approximation to $e$ in an SGFEM space $W\subset V$ that is richer than $X$, i.e., $W\supset X$. The quality of the resulting approximation is closely related to the quality of the Galerkin approximation $u_W\in W$ satisfying
\begin{align}\label{para_discrete_rich_prob}
\text{find }u_W\in W:\quad B(u_W,v) = F(v),\quad\text{for all }v \in W.
\end{align}

By letting $e_W = u_W - u_X$ we see that
\begin{align}\label{rich_error_approx}
B(e_W,v) = B(u_W,v) - B(u_X,v) = F(v) - B(u_X,v),\quad\text{for all } v\in W,
\end{align}
and thus $e_W\in W$ satisfying \eqref{rich_error_approx} estimates the true error $e\in V$. Clearly, since $e_W$ estimates $e$, SGFEM spaces $W$ that contain significantly improved approximations $u_W$ to $u$ (compared to $u_X$), also contain good estimates $e_W$ to $e$. To analyse the quality of the error estimate $||e_W||_B$, for a given choice of $W$, we require the following assumption.
\smallskip
\begin{assump}\label{Assump:saturation}\normalfont Let the functions $u$, $u_X$ and $u_W$ satisfy \eqref{para_weak_prob}, \eqref{para_discrete_prob} and \eqref{para_discrete_rich_prob} respectively. There exists a constant $\beta\in[0,1)$ (the saturation constant) such that
\begin{align}\label{saturation_bound}
||u-u_W||_B \leq\beta ||u-u_X||_B.
\end{align}
\end{assump}\noindent
We will also assume that $W:= X\oplus Y$ for some space $Y\subset V$ (the `detail' space) such that $X\cap Y = \{0\}$. Since computing $e_W\in W$ satisfying \eqref{rich_error_approx} is usually too expensive we instead exploit the decomposition of $W$ and solve:
\begin{align}\label{cheap_error_approx}
\text{find } e_Y\in Y:\quad B_0(e_Y,v) = F(v) - B(u_X,v),\quad\text{for all } v\in Y.
\end{align}
Notice the use of the parameter-free $B_0(\cdot,\cdot)$ bilinear form from \eqref{B0_bilinear_form} on the left-hand side of \eqref{cheap_error_approx}. To analyse the quality of the approximation $||e_Y||_{B_0}\approx ||e_W||_B$ we require the following result. Since $X$ and $Y$ are disjoint, and $B_0(\cdot,\cdot)$ induces a norm on the Hilbert space $V$ in \eqref{para_weak_prob}, there exists a constant $\gamma\in[0,1)$ such that
\begin{align}\label{B0_SCS}
|B_0(u,v)| \leq \gamma ||u||_{B_0}||v||_{B_0},\quad\text{for all }u\in X,\quad\text{for all } v \in Y,
\end{align}
see \cite[Theorem 5.4]{MR1885308}. Utilising \eqref{saturation_bound} and \eqref{B0_SCS} yields the following result \cite{CBS, MR3519560}.
\smallskip
\begin{thm}\label{Thm:stoch_finalbound} Let $u\in V=H_{0}^{1}(D) \otimes L_{\pi}^{2}(\Gamma)$ satisfy the variational problem \eqref{para_weak_prob} associated with the parametric diffusion problem \eqref{intro_PDE_problem1}--\eqref{intro_PDE_problem2} and let $u_X\in X$ satisfy \eqref{para_discrete_prob} for $X$ in \eqref{multilevel_X}. Choose $Y\subset V$ such that $X\cap Y = \{0\}$ and let $e_Y \in Y$ satisfy \eqref{cheap_error_approx}. If Assumption \ref{Assump:saturation} holds, as well as Assumptions \ref{Assump:a_bounds} and \ref{Assump:a0_bounds}, then $\eta:= ||e_Y||_{B_0}$ satisfies 
\begin{align}\label{stoch_finalbound}
\sqrt{\lambda} \, \eta \leq ||u-u_X||_{B}\leq \frac{\sqrt{\Lambda}}{\sqrt{1-\gamma^2}\sqrt{1-\beta^2}} \, \eta,
\end{align}
where $\lambda$ and $\Lambda$ are defined in \eqref{lambdas_def}, $\gamma\in[0,1)$ satisfies \eqref{B0_SCS}, and $\beta \in[0,1)$ satisfies \eqref{saturation_bound}.
\end{thm}

\smallskip
\noindent The quality of the error estimate $\eta\approx||e||_B$ depends on our choice of $Y$ because the constants $\gamma$ and $\beta$ in \eqref{stoch_finalbound} depend on $Y$. In the next section we describe a suitable structure for $Y$ when $X$ has the multilevel structure in \eqref{multilevel_X}.

\subsection{Choice of Detail Space $Y$}\label{Sec:a_post_error_para_diff}
In order to compute $\eta = ||e_Y||_{B_0}$ by solving \eqref{cheap_error_approx}, we need to choose the space $Y$. Note that in an adaptive SGFEM algorithm, $Y$ must vary with $X$, which is enriched at each step as we reduce $||u-u_X||_B$.  Suppose that $X$ has the form \eqref{multilevel_X}, where $J_{P}$ and the set of finite element spaces $\mathbf{H}_{1}$ are given. As stated in \cite[Remark 4.3]{MR3519560}, one possibility is to choose a second set of multi-indices $J_Q\subset J$ that satisfy $J_Q\cap J_P = \emptyset$ and construct
\begin{align}\label{Y_construction}
Y := \bigg(\bigoplus_{\mu\in J_P} H_2^{\mu}\otimes P^{\mu}\bigg)\oplus \bigg(\bigoplus_{\nu\in J_Q} H\otimes P^{\nu}\bigg),
\end{align}
where $H_2^{\mu}\subset H_0^1(D)$ are FEM spaces satisfying $H_1^{\mu}\cap H_2^{\mu}=\{0\}$ for all $\mu\in J_P$ and $H\subset H_0^1(D)$ is some other finite element space (to be defined later). Clearly, we have 
\begin{align}\label{Y_construction2}
Y:=Y_1 \oplus Y_2 :=\bigg(\bigoplus_{\mu\in J_P} Y_1^{\mu}\bigg) \oplus \bigg(\bigoplus_{\nu\in J_Q} Y_2^{\nu}\bigg),\qquad Y_1^\mu := H_2^{\mu}\otimes P^{\mu},\qquad Y_2^\nu := H\otimes P^{\nu},
\end{align}
which in turn leads to the following decomposition of $e_Y\in Y$, 
\begin{align*}
e_Y = e_{Y_1}+e_{Y_2} = \sum_{\mu\in J_P} e_{Y_1}^{\mu} + \sum_{\nu\in J_Q} e_{Y_2}^{\nu},\qquad e_{Y_1}^\mu\in Y_1^{\mu},\qquad e_{Y_2}^\nu\in Y_2^{\nu}.
\end{align*}

Since $B_0(\cdot,\cdot)$ is parameter-free and $J_P\cap J_Q = \emptyset$, then, as a consequence of the orthogonality property \eqref{mutual_orthogonality}, problem \eqref{cheap_error_approx} decouples into $\text{card}(J_P\cup J_Q) = \text{card}(J_P) + \text{card}(J_Q)$ smaller problems:
\begin{align}
\text{find }e_{Y_1}^{\mu}&\in Y_1^\mu:& B_0(e_{Y_1}^{\mu},v) &= F(v) - B(u_X,v),&\text{for all } v\in Y_1^\mu,& &\mu\in J_P,\label{spatial_error_solves}\\
\text{find }e_{Y_2}^{\nu}&\in Y_2^\nu:& B_0(e_{Y_2}^{\nu},v) &= F(v) - B(u_X,v),&\text{for all } v\in Y_2^\nu,& &\nu\in J_Q. \label{para_error_solves}
\end{align}
In addition, the error estimate $\eta$ in \eqref{stoch_finalbound} admits the decomposition
\begin{align}\label{eta}
\eta = ||e_Y||_{B_0} = \big(||e_{Y_1}||_{B_0}^2 + ||e_{Y_2}||_{B_0}^2\big)^{\frac{1}{2}} =  \bigg(\sum_{\mu\in J_P}||e_{Y_1}^{\mu}||_{B_0}^2 + \sum_{\nu\in J_Q}||e_{Y_2}^{\nu}||_{B_0}^2\bigg)^{\frac{1}{2}}.
\end{align}
For each $\mu\in J_P$ in \eqref{spatial_error_solves} we solve a problem of size $N_{Y_1}^\mu := \text{dim}(H_2^{\mu}\otimes P^{\mu}) = \text{dim}(H_2^{\mu})$ . For each $\nu\in J_Q$ in \eqref{para_error_solves}, we solve a problem of size $N_{Y_2}^\nu := \text{dim}(H\otimes P^{\nu}) = \text{dim}(H)$. We refer to $||e_{Y_1}||_{B_0}$ as the \emph{spatial error estimate}, and to $||e_{Y_2}||_{B_0}$ as the \emph{parametric error estimate}. For the adaptive algorithm in Section \ref{Sec:5}, it will be beneficial to define the set $\mathbf{H}_{2} = \left\{H_{2}^{\mu} \right\}_{\mu \in J_{P}}$ as well as the sets
\begin{align*}
\mathbf{N}_{Y_1}=\{N_{Y_1}^\mu\}_{\mu\in J_P},\qquad \mathbf{N}_{Y_2}=\{N_{Y_2}^\nu\}_{\nu\in J_Q}.
\end{align*}
The quality of the error estimate $\eta$ depends on our choice of $J_Q$ and $\mathbf{H}_{2}$ as well as the finite element space $H$ appearing in the definition of $Y_{2}$, since they affect the constants $\gamma$ and $\beta$ appearing in \eqref{stoch_finalbound}. The error bound is sharp when $\beta$ and $\gamma$ are close to zero.
 
If Assumption \ref{Assump:a0_bounds} holds, then $H_0^1(D)$ is a Hilbert space with respect to the inner product 
\begin{align*}
\langle a_0u,v\rangle = \int_Da_0(\mathbf{x})\nabla u(\mathbf{x})\cdot\nabla v(\mathbf{x})\ d\mathbf{x},\qquad u,v\in H_0^1(D).
\end{align*}
Furthermore, since $H_1^{\mu}\cap H_2^{\mu} = \{0\}$ for all $\mu\in J_P$, there exists a constant $\gamma^{\mu}\in [0,1)$ such that
\begin{align}\label{a0_SCS}
|\langle a_0u,v\rangle| \leq \gamma^{\mu} \langle a_0u,u\rangle^{1/2}\langle a_0v,v\rangle^{1/2},\quad\text{for all } u \in H_1^{\mu},\quad\text{for all } v\in H_2^{\mu},
\end{align}
for all $\mu\in J_P$ (again, see \cite[Theorem 5.4]{MR1885308}).  We denote the smallest such constant (known as the CBS constant) by $\gamma_{\textrm{min}}^{\mu}$. Note that this constant only depends on the chosen finite element spaces $H_1^{\mu}$ and $H_2^{\mu}$ and is known explicitly in many cases, see \cite{CBS}. It is then straightforward to prove, using the mutual orthogonality of the sets $\left\{\psi_{\mu}(\mathbf{y} )\right\}_{\mu \in J_{P}}$ and $\left\{\psi_{\nu}(\mathbf{y} )\right\}_{\nu \in J_{Q}}$ and the definition of $B_{0}(\cdot, \cdot)$ that with $Y$ chosen as in \eqref{Y_construction},  the bound \eqref{B0_SCS} holds with 
\begin{align}\label{gamma_result}
\gamma:=\textrm{max}_{\mu \in J_{P}}\left\{\gamma_{\textrm{min}}^{\mu} \right\}.
\end{align}
See also \cite[Remark 4.3]{MR3519560}. 
\smallskip
\begin{rem}
Since $H$ in \eqref{Y_construction} does not depend on $\nu\in J_Q$, the matrix that characterises the linear systems associated with \eqref{para_error_solves} is the same for all $\nu\in J_Q$. Only the right-hand side changes. Consequently, we can vectorise the system solves associated with \eqref{para_error_solves} over the multi-indices $J_Q$.
\end{rem}
\begin{rem}\label{Rem:CBSrem2} For two FEM spaces $H_1^\mu$ and $H_2^\mu$, there often exists a sharp upper bound for the associated CBS constant $\gamma_{\min}^\mu$ that is independent of the mesh level number $\ell^\mu$, see \cite{CBS}.
\end{rem}

\subsection{The Spatial Error Estimator}\label{Sec:spatial_apost}
\begin{figure}[t!]
\centering
\subfigure{
\begin{minipage}[c]{0.24\linewidth}
\centering
\begin{tikzpicture}
\fill[lightgray] (0,0) rectangle (2.5,2.5);
\draw[step=1.25cm,darkgray,line width = 0.6mm] (0,0) grid (2.5,2.5);
\filldraw[fill=white, draw = black, very thick](0   ,0   ) circle (0.12cm);
\filldraw[fill=white, draw = black, very thick](1.25,0   ) circle (0.12cm);
\filldraw[fill=white, draw = black, very thick](2.5 ,0   ) circle (0.12cm);
\filldraw[fill=white, draw = black, very thick](0   ,1.25) circle (0.12cm);
\filldraw[fill=white, draw = black, very thick](1.25,1.25) circle (0.12cm);
\filldraw[fill=white, draw = black, very thick](2.5 ,1.25) circle (0.12cm);
\filldraw[fill=white, draw = black, very thick](2.5 ,0   ) circle (0.12cm);
\filldraw[fill=white, draw = black, very thick](0   ,2.5 ) circle (0.12cm);
\filldraw[fill=white, draw = black, very thick](1.25,2.5 ) circle (0.12cm);
\filldraw[fill=white, draw = black, very thick](2.5 ,2.5 ) circle (0.12cm);
\end{tikzpicture}
\end{minipage}
}
\subfigure{
\begin{minipage}[c]{0.24\linewidth}
\centering
\begin{tikzpicture}
\fill[lightgray] (0,0) rectangle (2.5,2.5);
\draw[step=1.25cm,darkgray,line width = 0.6mm] (0,0) grid (2.5,2.5);
\filldraw[fill=white, draw = black, very thick](0.625,0    ) circle (0.12cm);
\filldraw[fill=white, draw = black, very thick](1.875,0    ) circle (0.12cm);
\filldraw[fill=white, draw = black, very thick](0    ,0.625) circle (0.12cm);
\filldraw[fill=white, draw = black, very thick](0.625,0.625) circle (0.12cm);
\filldraw[fill=white, draw = black, very thick](1.25 ,0.625) circle (0.12cm);
\filldraw[fill=white, draw = black, very thick](1.875,0.625) circle (0.12cm);
\filldraw[fill=white, draw = black, very thick](2.5  ,0.625) circle (0.12cm);
\filldraw[fill=white, draw = black, very thick](0.625,1.25 ) circle (0.12cm);
\filldraw[fill=white, draw = black, very thick](1.875,1.25 ) circle (0.12cm);
\filldraw[fill=white, draw = black, very thick](0    ,0.625) circle (0.12cm);
\filldraw[fill=white, draw = black, very thick](0    ,1.875) circle (0.12cm);
\filldraw[fill=white, draw = black, very thick](0.625,1.875) circle (0.12cm);
\filldraw[fill=white, draw = black, very thick](1.25 ,1.875) circle (0.12cm);
\filldraw[fill=white, draw = black, very thick](1.875,1.875) circle (0.12cm);
\filldraw[fill=white, draw = black, very thick](2.5  ,1.875) circle (0.12cm);
\filldraw[fill=white, draw = black, very thick](0.625,2.5  ) circle (0.12cm);
\filldraw[fill=white, draw = black, very thick](1.875,2.5  ) circle (0.12cm);
\end{tikzpicture}
\end{minipage}
}
\caption{The nodes associated with $H_1^\mu$ (left) and $H_2^\mu$ (right), when $H_1^\mu$ is chosen to be a $\mathbb{Q}_1$ space and $H_2^\mu$ is chosen to be a `broken' $\mathbb{Q}_2$ space associated with the same mesh $\mathcal{T}_{\ell^\mu}$ as $H_1^\mu$.}\label{Fig:H2detail}
\end{figure}
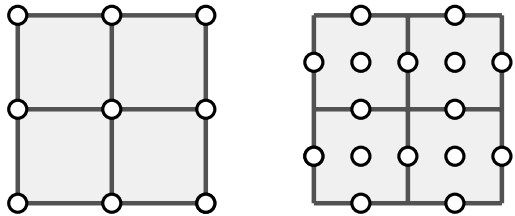

We now briefly discuss possible choices of the FEM spaces $\mathbf{H}_2 = \{H_2^\mu\}_{\mu\in J_P}$ that define the tensor spaces $\mathbf{Y}_1 := \{Y_1^{\mu}\}_{\mu\in J_P}$ in \eqref{Y_construction2}. Recall that each FEM space $H_1^\mu$ is associated with a mesh $\mathcal{T}_{\ell^\mu}= \mathcal{T}_{i}$ for some $i\in\mathbb{N}_0$. One option is to construct a basis for $H_2^\mu$ with respect to the same mesh $\mathcal{T}_{\ell^\mu}$ but using polynomials of a higher degree. In order to ensure that $H_1^\mu\cap H_2^\mu = \{0\}$, we exclude basis functions associated with nodes associated with $H_{1}^{\mu}$. For example, if the spaces $\mathbf{H}_1$ are $\mathbb{Q}_1$ FEM spaces, we may choose the spaces $\mathbf{H}_2$ to be `broken' $\mathbb{Q}_2$ FEM spaces (see Figure \ref{Fig:H2detail}). Another option is to use polynomials of the same degree, but introduce basis functions associated with the new nodes that would be introduced by performing the mesh refinement $\mathcal{T}_{\ell^\mu}\to \mathcal{T}_{i+1}$ (i.e., by increasing the level number by one).

\subsection{The Parametric Error Estimator}\label{Sec:parametric_apost}
It remains to explain how to choose the multi-indices $J_Q$ and the space $H \subset H_0^1(D)$ that define the tensor spaces $\mathbf{Y}_2 := \{Y_2^{\mu}\}_{\mu\in J_Q}$ in \eqref{Y_construction2}. It was proven in \cite{MR3519560} that $||e_{Y_2}^\nu||_{B_0} = 0$ for considerably many multi-indices $\nu\in J\backslash J_P$. In order to avoid unnecessary computations, it is essential that we first identify the set of multi-indices $J^*\subset J$ that result in non-zero contributions. Indeed, this set is given by
\begin{align*}
J^* = \big\{\mu\in J\backslash J_P;\ \mu=\nu \pm \epsilon^m\ \forall\ \nu\in J_P,\ \forall\ m\in\mathbb{N} \big\},
\end{align*}
where $\epsilon^m : = (\epsilon_1^m,\epsilon_2^m,\dots)$ is the Kronecker delta sequence such that $\epsilon_j^m = \delta_{mj}$ for all $j\in\mathbb{N}$. Since $J^*$ is an infinite set, we need to choose a finite subset $J_Q\subset J^*$. We call $J^*$ the set of `neighbouring indices' to $J_P$ and choose
\begin{align}\label{JQ_setup}
J_Q = \big\{\nu\in J^*;\ \max\{\text{supp}(\nu)\}\leq M + \Delta_M \big\},
\end{align}
where $\Delta_M\in\mathbb{N}$ is  the number of additional parameters we wish to activate.


We now turn our attention to $H\subset H_0^1(D)$. Recall that $W = X\oplus Y$ in \eqref{para_discrete_rich_prob}.  The space $Y$ (and hence $Y_{2}^{\nu}= H \otimes P^{\nu}$) should be chosen so that $W$ contains functions that would result in an improved approximation $u_W\in W$ to $u$. We clearly want to choose $Y$ so that we have an accurate energy error estimate $\eta$ for the current approximation $u_{X}$. However, since we want to perform adaptivity, the functions in $Y$ serve as candidates to be added to $X$ at the next approximation step.  Since $X$ may be augmented with $H\otimes P^\nu$ for some $\nu\in J_Q$, we should choose $H$ such that the structure of $Y$ in \eqref{Y_construction} is maintained and the error estimator is straightforward to compute at each step. For this reason, we choose $H = H_1^{\bar\mu}$ for some $\bar\mu\in J_P$. That is, we choose $H$ to be one of the FEM spaces already used in the construction of $X$.

When choosing $\bar\mu \in J_P$ we must consider the fact that through our choice of $Y$ in \eqref{Y_construction}, $\beta$ in \eqref{stoch_finalbound} depends on $\bar\mu$. We have to balance the accuracy of the estimate $\eta$ against the cost to compute it. If we choose $\bar\mu$ such that $\ell^{\bar\mu}=\max_{\mu\in J_P}\bm\ell$ (i.e., choose the richest FEM space used so far), then $\textrm{dim}(X)$ will grow too quickly when we augment $X$ with functions in $\mathbf{Y}_2$. Similarly, if $\ell^{\bar \mu}= \min_{\mu\in J_P}\bm\ell$, the error reduction may be negligible if $X$ is augmented with functions from $\mathbf{Y}_2$. To strike a balance, we will choose $\bar\mu$ to correspond to the FEM space $H_{1}^{\mu}$ with the smallest mesh level number $\ell^\mu$ such that the number of spaces with level number $\ell^{\mu}$ or less is greater than or equal to $\lceil\frac{1}{2}\text{card}(J_P)\rceil$. We denote this choice by $\bar\mu = \argavg_{\mu\in J_P}\bm\ell$.
\smallskip
\begin{examp}\label{Examp:mubar}\normalfont Suppose $\textrm{card}(J_P)=5$ and $\bm\ell=\{2,3,3,2,1\}$, then $\ell^{\bar\mu}= 2$. Similarly, if $\textrm{card}(J_P)=3$ and $\bm\ell = \{4,3,2\}$, then $\ell^{\bar\mu}= 3$. 
\end{examp}

\smallskip 
\noindent The choice $\bar\mu=\argavg_{\mu\in J_P}\bm\ell$ ensures that the dimensions of the spaces in $\mathbf{Y}_2$ are always modest in comparison to those of the spaces in $\mathbf{X}=\{X^\mu\}_{\mu\in J_P}$ in \eqref{multilevel_X}.

\section{Adaptive Multilevel SGFEM} \label{Sec:5}
Suppose that $X$ and $Y$ in \eqref{multilevel_X} and \eqref{Y_construction} have been chosen (and so the sets of multi-indices $J_P,J_Q\subset J$ have also been chosen) and that the corresponding approximations $u_X\in X$ and $e_Y\in Y$ satisfying \eqref{para_discrete_prob} and \eqref{cheap_error_approx} have been computed. If $\eta =||e_{Y}||_{B_{0}}$ is too large, we want to augment $X$ with some of the functions in $Y$ and compute a (hopefully) improved approximation to $u\in V$ satisfying \eqref{para_weak_prob}. Of course, we could augment $X$ with the full space $Y$ to ensure it is sufficiently rich. However, we must also ensure that the total number of additional degrees of freedom (DOFs) introduced is balanced against the reduction in the energy error that is achieved. We should only augment $X$ with functions that result in significant error reductions. Below, we demonstrate that using the sets of component estimates
\begin{align}\label{err_component_sets}
\mathbf{E}_{Y_1}:=\{||e_{Y_1}^{\mu}||_{B_0}\}_{\mu\in J_P},\qquad \mathbf{E}_{Y_2}:=\{||e_{Y_2}^{\mu}||_{B_0}\}_{\mu\in J_Q},
\end{align}
(which are computed to determine $\eta$), we can estimate the error reduction that would be achieved by performing certain enrichment strategies at the next approximation step.

\subsection{Estimated Error Reductions}\label{Sec:enrich_strat}
Consider the discrete problems:
\begin{align}
\text{find } u_{W_1} \in W_1:\qquad B(u_{W_1},v)&= F(v),\qquad \text{for all } v \in W_1,\label{W1_enrich}\\
\text{find } u_{W_2} \in W_2:\qquad B(u_{W_2},v)&= F(v),\qquad \text{for all } v \in W_2,\label{W2_enrich}
\end{align}
where $W_1$ and $W_2$ are `enhanced' SGFEM approximation spaces given by
\begin{equation}\label{enrichment_strategies}
\begin{aligned}
W_1 &:= X\oplus Y_{W_1} := X\oplus \bigg(\bigoplus_{\mu\in \bar J_P} Y_1^\mu \bigg),&
& \bar J_P \subset J_P,\\
W_2 &:= X\oplus Y_{W_2} := X\oplus \bigg(\bigoplus_{\nu\in \bar J_Q} Y_2^\nu \bigg),&
& \bar J_Q \subset J_Q.
\end{aligned}
\end{equation}
That is, $u_{W_1}$ and $u_{W_2}$ are SGFEM approximations to $u \in V$ computed in $W_1$ and $W_2$, respectively. Note that if $\bar J_P = J_P$ then $Y_{W_{1}}=Y_{1}$ and if  $\bar J_Q = J_Q$ then $Y_{W_{2}}=Y_{2}$. However, we want to consider enrichment strategies associated with only important subsets of the multi-indices. The space $W_{1}$ corresponds to refining the finite element meshes associated with a subset of the multi-indices $\mu \in J_{P}$ used in the definition of $X$, whereas $W_{2}$ corresponds to adding new basis polynomials on the parameter domain. We want to estimate the potential pay-offs of these two strategies.

Let $e_{W_1} = u - u_{W_1}$ denote the error corresponding to the enhanced approximation $u_{W_1}$. Due to the orthogonality of $e_{W_1}$ with functions in $W_1$ ($(u_{W_1} - u_X)\in W_1$ in particular) with respect to $B(\cdot,\cdot)$ (Galerkin-orthogonality), and the symmetry of $B(\cdot,\cdot)$, we find that
\begin{align*}
||e_{W_1}||_B^2 = ||u-u_X||_B^2 - ||u_{W_1} - u_X||_B^2.
\end{align*}
Hence, $||u_{W_1} - u_X||_B^2$ characterises the reduction in $||u-u_X||_B^2$ (the square of the energy error) that would be achieved by augmenting $X$ with $Y_{W_1}$, for a suitably chosen set $\bar J_P \subset J_P$, and computing an enhanced approximation $u_{W_1}\in W_1$ satisfying \eqref{W1_enrich}. Similarly, $||u_{W_2} - u_X||_B^2$ characterises the reduction in $||u-u_X||_B^2$ that would be achieved by augmenting $X$ with $Y_{W_2}$ for a suitably chosen set $\bar J_Q\subset J_Q$ and computing $u_{W_2}\in W_2$ satisfying \eqref{W2_enrich}. The following result provides estimates for these quantities. This is a simple extension of a result proved in \cite{MR3177362,MR3519560}; the proof is very similar.

\smallskip
\begin{thm}\label{Thm:est_error_reduc}\normalfont Let $u_X\in X$ be the SGFEM approximation satisfying \eqref{para_discrete_prob} and let $u_{W_1}\in W_1$ and $u_{W_2}\in W_2$ satisfy problems \eqref{W1_enrich} and \eqref{W2_enrich}. Define the quantities
\begin{align*}
\zeta_{W_1} := \sum_{\mu\in \bar J_P}||e_{Y_1}^\mu||_{B_0}^2,\qquad \zeta_{W_2} := \sum_{\nu\in \bar J_Q}||e_{Y_2}^\nu||_{B_0}^2,
\end{align*}
for some $\bar J_P\subset J_P$ and $\bar J_Q\subset J_Q$. Then the following estimates hold:
\begin{align}
&\lambda\zeta_{W_1}\leq ||u_{W_1} - u_X||_{B}^2\leq \frac{\Lambda}{1-\gamma^2}\zeta_{W_1},\label{error_reduc_1}\\
&\lambda\zeta_{W_2}\leq ||u_{W_2} - u_X||_{B}^2\leq \Lambda\zeta_{W_2},\label{error_reduc_2}
\end{align}
where $\lambda$ and $\Lambda$ are the constants in \eqref{lambdas_def}, and $\gamma\in[0,1)$ is the constant satisfying \eqref{gamma_result}.
\end{thm}

Given two sets of multi-indices $\bar J_P$ and $\bar J_Q$, we now determine an appropriate enrichment strategy for $X$ by considering the bounds \eqref{error_reduc_1}--\eqref{error_reduc_2}. One option would be to perform the enrichment strategy that corresponds to $\max\{\zeta_{W_1},\zeta_{W_2}\}$. Whilst this may lead to a large reduction of $||u-u_X||_B^2$ (and hence of $||u-u_X||_B$), it doesn't take into account the computational cost incurred. We want to construct sequences of SGFEM spaces $X$ for which the energy error converges to zero at the best possible rate with respect to $N_{\textrm{dof}}=\textrm{dim}(X)$ for the chosen set of finite element spaces.  Hence, the number of DOFs should be taken into account. Recall the definitions
\begin{align}\label{error_component_dimensions}
N_{Y_1}^\mu := \text{dim}(Y_1^\mu),\ \mu\in J_P,\qquad N_{Y_2}^\nu := \text{dim}(Y_2^\nu),\ \nu\in J_Q.
\end{align}
The number of additional DOFs (compared to the current space $X$) associated with the spaces $W_1$ and $W_2$ in \eqref{enrichment_strategies} is given by 
\begin{align*}
N_{W_1} := \sum_{\mu\in \bar J_P} N_{Y_1}^\mu,\qquad N_{W_2} := \sum_{\nu\in \bar J_Q} N_{Y_2}^\nu,
\end{align*}
respectively. Due to Theorem \ref{Thm:est_error_reduc}, the ratios
\begin{align}\label{final_ratios}
R_{ W_1}:=\frac{\zeta_{W_1}}{N_{W_1}},\qquad
R_{ W_2}:=\frac{\zeta_{W_2}}{N_{W_2}},
\end{align}
provide approximations to $||u_{W_1}-u_X||_B^2/N_{W_1}$ and $||u_{W_2}-u_X||_B^2/N_{W_2}$, respectively. Once we have chosen $\bar J_P$ and $\bar J_Q$, we augment $X$ with the space $Y_{W_1}$ or $Y_{W_2}$, that corresponds to $\max\{R_{W_1},R_{W_2}\}$. In the next section we propose an adaptive multilevel SGFEM algorithm for the numerical solution of \eqref{intro_PDE_problem1}--\eqref{intro_PDE_problem2} as well as two methods for the selection of the sets of multi-indices $\bar J_P$ and $\bar J_Q$.

\begin{algorithm}[t!]
    \caption{Adaptive multilevel SGFEM}\label{Alg:1}
	\medskip
    \SetKwInOut{Input}{Input}
    \SetKwInOut{Output}{Output}
    \Input{Problem data $a(\mathbf{x},\mathbf{y})$, $f(\mathbf{x})$; initial index set $J_P^0$ and mesh level numbers $\bm\ell^0$;
    energy error tolerance $\epsilon$.}
    \Output{Final SGFEM approximation $u_X^K$ and energy error estimate $\eta^K$.}
    \medskip
    Choose \texttt{version} (1 or 2)\\
	\For{$k=0,1,2,\dots$}
	{
	$u_X^k\leftarrow$ \texttt{SOLVE}$\big[a,f,J_P^k,\bm\ell^k\big]$\\
	$J_Q^k\leftarrow$ \texttt{PARAMETRIC\_INDICES}$\big[J_P^k\big]$\hfill see:\hspace{3pt} \eqref{JQ_setup}\\
	$\mathbf{E}_{Y_1}^k\leftarrow$ \texttt{COMPONENT\_SPATIAL\_ERRORS}$\big[u_X^k,J_P^k,\bm\ell^k\big]$\hfill \eqref{err_component_sets}\\
		$\mathbf{E}_{Y_2}^k\leftarrow$ \texttt{COMPONENT\_PARAMETRIC\_ERRORS}$\big[u_X^k,J_Q^k,\bm\ell^k\big]$\\
	$\eta^k = \big[\sum_{\mu\in J_P^k}||e_{Y_1}^{\mu,k}||_{B_0}^2+\sum_{\nu\in J_Q^k}||e_{Y_2}^{\nu,k}||_{B_0}^2\big]^{\frac{1}{2}}$
	\hfill \eqref{eta}\\
	\eIf{$\eta^k<\epsilon$}
	{\textbf{return} $u_X^k,\eta^k$}
	{$[\texttt{refinement\_type},\bar J^k]\leftarrow \texttt{ENRICHMENT\_INDICES}\big[\texttt{version},\mathbf{E}_{Y_1}^k,\mathbf{E}_{Y_2}^k,J_P^k,J_Q^k\big]$\\
	\eIf{\emph{\texttt{refinement\_type}}$\ =\ $\emph{\texttt{spatial}}}
	{$J_P^{k+1} = J_P^k$\\
	$\bm\ell^{k+1} = \big\{\ell_k^{\mu+};\ \mu\in \bar J^k\big\}\cup\big\{\ell_k^\mu;\ \mu \in J_P^k\backslash \bar J^k\big\}$ \hfill \eqref{h_mu_plus}
	}
	{$J_P^{k+1} = J_P^k\cup \bar J^k$\\
	$\bm\ell^{k+1} = \bm\ell^k \cup \big\{\ell_k^{\bar\mu};\ \nu\in\bar J^k\big\}$
	}     
     }
    
	}
\end{algorithm}

\subsection{An Adaptive Algorithm}
Using the a posteriori error estimation strategy discussed in Section \ref{Sec:a_post_error_para_diff}, and the estimated error reductions described in Section \ref{Sec:enrich_strat}, we now propose an adaptive algorithm that generates a sequence of multilevel SGFEM spaces
\begin{align*}
X^0 \subset X^1 \cdots \subset X^k \cdots \subset X^K \subset V,
\end{align*}
and terminates at step $k = K$ when the SGFEM approximation $u_X^K\in X^K$ to $u$ satisfies a prescribed error tolerance $\epsilon$. We start by selecting an initial low-dimensional SGFEM space $X^0$ of the form \eqref{multilevel_X} and compute an initial approximation $u_X^0\in X^0$ to $u\in V$ satisfying \eqref{para_discrete_prob}. Assuming that the polynomial degree of the FEM approximation on $D$ has been fixed,  we only need to supply an initial set of multi-indices $J_P^0$, as well as a set of mesh level numbers $\bm\ell^0 = \{\ell_0^\mu\}_{\mu\in J_P^0}$.
We then consider two enrichment strategies.
The first option is to refine certain meshes associated with the spaces $\mathbf{H}_1^0$ and produce a new set $\bm\ell^1$. If $\ell_0^\mu = i$ for some $\mu\in J_P$, and we want to perform a refinement, we set $\ell_1^\mu = i+1$ or equivalently replace $\mathcal{T}_{\ell_0^\mu}$ with the next mesh in the sequence $\bm{\mathcal{T}}$ in \eqref{sequence_of_meshes}. In our adaptive algorithm we write
\begin{align}\label{h_mu_plus}
\ell_0^\mu \to \ell_0^{\mu+} =: \ell_1^\mu.
\end{align}
The second option is to add multi-indices to $J_P^0$ to give a new set $J_{P}^{1}$. In this case, we must also update $\bm\ell^0$ with new mesh parameters to maintain the relationship $\text{card}(J_P)=\text{card}(\bm\ell)$. Specifically, we add a copy of $\ell_0^{\bar\mu}$ to $\bm\ell^0$, for every multi-index added to $J_P^0$ (see Section \ref{Sec:parametric_apost} for the definition of $\bar\mu$). Once $J_P^1$ and $\bm\ell^1$ are defined, and $u_X^1\in X^1$ is computed, the process is repeated.

\begin{algorithm}[ht!]
    \caption{\texttt{ENRICHMENT\_INDICES} versions 1 and 2}\label{Alg:2}
    \medskip
    \SetKwInOut{Input}{Input}
    \SetKwInOut{Output}{Output}
    \Input{\texttt{version}; $\mathbf{E}_{Y_1}^k$; $\mathbf{E}_{Y_2}^k$; $J_P^k$; $J_Q^k$.}
    \Output{\texttt{refinement\_type}, $\bar J^k$.}
    \medskip
    $\delta_{Y_1}^k = \max_{\mu\in J_P^k} \mathbf{R}_{Y_1}^k$, $\delta_{Y_2}^k = \max_{\nu\in J_Q^k} \mathbf{R}_{Y_2}^k$\\
    \eIf{$\delta_{Y_1}^k > \delta_{Y_2}^k$}
    {
    $\bar J_Q^k = \{\nu\in J_Q^k;\ R_{Y_2}^{\nu,k} = \delta_{Y_2}^k\}$\\
    \eIf{$\texttt{\emph{version}} = 1$}
    {$\bar J_P^k = \{\mu\in J_P^k;\ R_{Y_1}^{\mu,k} > \delta_{Y_2}^k\}$}
    {$\bar J_P^k\leftarrow\texttt{MARK}[\mathbf{E}_{Y_1}^k,\mathbf{N}_{Y_1}^k,\delta_{Y_2}^k]$}
    }
    {
$\bar J_P^k = \{\mu\in J_P^k;\ R_{Y_1}^{\mu,k} = \delta_{Y_1}^k\}$\\
    \eIf{$\texttt{\emph{version}} = 1$}
    {$\bar J_Q^k = \{\nu\in J_Q^k;\ R_{Y_2}^{\nu,k} > \delta_{Y_1}^k\}$}
    {$\bar J_Q^k\leftarrow\texttt{MARK}[\mathbf{E}_{Y_2}^k,\mathbf{N}_{Y_2}^k,\delta_{Y_1}^k]$}
}
	\eIf{$R_{W_1}^k > R_{W_2}^k$}
	{\texttt{refinement\_type} = \texttt{spatial},\ $\bar J^k = \bar J_P^k$}
	{\texttt{refinement\_type} = \texttt{parametric},\ $\bar J^k = \bar J_Q^k$}
\textbf{return} [\texttt{refinement\_type}, $\bar J^k$]
\end{algorithm}

The general process is outlined in Algorithm \ref{Alg:1}. At a given step $k$: 
\begin{itemize}
\item \texttt{SOLVE} computes an SGFEM approximation $u_X\in X$ to $u\in V$ satisfying \eqref{para_discrete_prob}.
\item \texttt{PARAMETRIC\_INDICES} uses \eqref{JQ_setup} to determine a subset $J_Q$ of the neighbouring indices to $J_P$ for a prescribed choice of $\Delta_M$.
\item \texttt{COMPONENT\_SPATIAL\_ERRORS} and \texttt{COMPONENT\_PARAMETRIC\_ERRORS} compute the sets of error estimates $\mathbf{E}_{Y_1}$ and $\mathbf{E}_{Y_2}$ in \eqref{err_component_sets}, respectively, by solving \eqref{spatial_error_solves} and \eqref{para_error_solves}.
\item \texttt{ENRICHMENT\_INDICES} analyses the sets $\mathbf{E}_{Y_1}$ and $\mathbf{E}_{Y_2}$ in conjunction with the formulae in \eqref{final_ratios} to determine how to enrich the current SGFEM space $X$.
\end{itemize} 
A key part of \texttt{ENRICHMENT\_INDICES} is the determination of suitable sets $\bar J_P\subset J_P$ and $\bar J_Q\subset J_Q$, which we describe in the next section. Algorithm \ref{Alg:1} subsequently performs either a spatial or parametric refinement associated with the set of multi-indices $\bar J := \bar J_P$ or $\bar J := \bar J_Q$, respectively.

\subsection{Selection of the Enrichment Multi-indices} We introduce two versions of the module \texttt{ENRICHMENT\_INDICES}, which are outlined in Algorithm \ref{Alg:2}. To begin, define the sets
\begin{align*}
\mathbf{R}_{Y_1} := \big\{ R_{Y_1}^{\mu} \big\}_{\mu\in J_P} := \bigg\{\frac{||e_{Y_1}^{\mu}||_{B_0}^2}{N_{Y_1}^{\mu}}\bigg\}_{\mu\in J_P},\quad
\mathbf{R}_{Y_2} := \big\{ R_{Y_2}^{\nu} \big\}_{\nu\in J_Q} := \bigg\{\frac{||e_{Y_2}^{\nu}||_{B_0}^2}{N_{Y_2}^{\nu}}\bigg\}_{\nu\in J_Q},
\end{align*}
of estimated error reduction ratios and consider the quantities
\begin{align*}
\delta_{Y_1}:= \max_{\mu\in J_P} \mathbf{R}_{Y_1},\qquad \delta_{Y_2}:= \max_{\nu\in J_Q}\mathbf{R}_{Y_2}.
\end{align*}

Version 1 of Algorithm 2 is simple. If $\delta_{Y_1} > \delta_{Y_2}$, we define $\bar J_P$ to be the set of multi-indices $\mu\in J_P$ such that $R_{Y_1}^{\mu} > \delta_{Y_2}$ and we define $\bar J_Q$ to be the set of multi-indices $\nu\in J_Q$ such that $R_{Y_2}^{\nu} = \delta_{Y_2}$. Similarly, if $\delta_{Y_2} > \delta_{Y_1}$, we define $\bar J_Q$ to be the set of multi-indices in $J_Q$ such that $R_{Y_2}^{\nu} > \delta_{Y_1}$  and $\bar J_P$ is the set of multi-indices in $J_P$ such that $R_{Y_1}^{\mu} =\delta_{Y_1}$. The refinement type is then determined by computing $R_{W_1}$ and $R_{W_2}$ in \eqref{final_ratios}. If $R_{W_1}>R_{W_2}$ we perform \emph{spatial} refinement and set $\bar J=\bar J_{P}$. Otherwise, we enrich the \emph{parametric} part, and set $\bar J=\bar J_{Q}$.

Version 2 is similar. However,  if $\delta_{Y_1} > \delta_{Y_2}$, we choose $\bar J_P$ to be the largest subset of $J_{P}$ such that $R_{W_1}>\delta_{Y_2}$ (recall $R_{W_1}$ depends on $\bar J_P$). Similarly, if $\delta_{Y_2} > \delta_{Y_1}$, we choose $\bar J_Q$ to be the largest subset of $J_{Q}$ such that $R_{W_2}>\delta_{Y_1}$. As before, the refinement type chosen is the one associated with $\max\{R_{W_1},R_{W_2}\}$. Version 2 is reminiscent of a D\"{o}rfler marking strategy \cite{MR1393904} and so the module that generates $\bar J_P$ (if $\delta_{Y_{1}} > \delta_{Y_{2}}$ ) and $\bar J_Q$ (if $\delta_{Y_{2}} > \delta_{Y_{1}}$) is called \texttt{MARK}.

\smallskip
\begin{rem} A key feature of both versions of \texttt{ENRICHMENT\_INDICES} is that no marking or tuning parameters are required. The user only needs to choose $\Delta_{M}$ in the definition of $J_{Q}$ in \eqref{JQ_setup}. This fixes an upper bound on the number of new parameters $y_{m}$ that may be activated. 
\end{rem}

\section{Numerical Experiments}\label{Sec:num_results}
We now investigate the performance of Algorithms \ref{Alg:1} and \ref{Alg:2} in computing approximate solutions to \eqref{intro_PDE_problem1}--\eqref{intro_PDE_problem2}. First, we describe four test problems. These differ, in particular, in the choice of $a(\mathbf{x},\mathbf{y})$, and give rise to sequences of coefficients $\left\{\|a_{m}\|_{\infty}\right\}_{m=1}^{\infty}$ that decay at different rates. Recall,  $y_{m} \in \Gamma_{m} =[-1,1]$ is the image of a uniform random variable and $\pi_{m}(y_{m})$ is the associated probability measure, for $m \in \mathbb{N}$. 

\subsubsection*{Test Problem 1 (TP.1)}
First, we consider a problem from \cite{MR3177362,CBS}. Let $f(\mathbf{x})=\frac{1}{8}(2-x_1^2-x_2^2)$ for $\mathbf{x} = (x_1,x_2)^{\top}\in D:=[-1,1]^2$ and assume that 
\begin{align}\label{TP1_a}
a(\mathbf{x},\mathbf{y}) = 1 + \sigma\sqrt{3}\sum_{m=1}^\infty\sqrt{\lambda_m}\phi_m(\mathbf{x})y_m,
\end{align}
where $(\lambda_m,\phi_m)$ are the eigenpairs of the operator associated with the covariance function
\begin{align*}
C[a](\mathbf{x},\mathbf{x}') = \exp\bigg(-\frac{|x_1-x_1'|}{l_1}-\frac{|x_2-x_2'|}{l_2}\bigg),\qquad \mathbf{x},\mathbf{x}'\in D.
\end{align*}
As in \cite{CBS} we choose $\sigma =0.15$ (the standard deviation) and $l_1=l_2=2$ (the correlation lengths). It can be shown that asymptotically (as $m \to \infty$), $\lambda_m$ is  $\mathcal{O}(m^{-2})$, see \cite{MR3308418}. 

\subsubsection*{Test Problem 2 (TP.2)} Next, we consider a problem from \cite{MR3154028,MR3519560}. Let $f(\mathbf{x})=1$ for $\mathbf{x} = (x_1,x_2)^{\top}\in D:=[0,1]^2$ and assume that
\begin{align*}
a(\mathbf{x},\mathbf{y}) = 1 + \sum_{m=1}^\infty \alpha_m\cos(2\pi\beta_m^1x_1)\cos(2\pi\beta_m^2x_2)y_m,
\end{align*}
where $\beta_m^1 = m-k_m(k_m+1)/2$, $\beta_m^2 = k_m - \beta_m^1$ and $k_m = \lfloor -1/2+(1/4 + 2m)^{1/2} \rfloor$ for $m\in\mathbb{N}$. In this test problem, we select the amplitude coefficients $\alpha_m = 0.547m^{-2}$.
\subsubsection*{Test Problem 3 (TP.3)} This is the same as TP.2  but we now choose $\alpha_m = 0.832m^{-4}$, so that the terms in the expansion of $a(\mathbf{x},\mathbf{y})$ decay more quickly.

\subsubsection*{Test Problem 4 (TP.4)} Finally, we consider a problem from \cite{MR3308418}. Let $f(\mathbf{x})$ and $D$ be as in TP.2 and assume that
\begin{align}\label{TP4_a}
a(\mathbf{x},\mathbf{y}) = 2 + \sqrt{3}\sum_{i=0}^\infty\sum_{j=0}^\infty \sqrt{\nu_{ij}}\phi_{ij}(\mathbf{x})y_{ij} ,\qquad y_{ij}\in [-1,1]
\end{align}
where $\phi_{00}= 1$, $\nu_{00} = \frac{1}{4}$ and
\begin{align*}
\phi_{ij} = 2\cos(i\pi x_1)\cos(j\pi x_2),\quad \nu_{ij} = \frac{1}{4}\exp(-\pi(i^2+j^2)l^{-2}).
\end{align*}
We choose the correlation length $l=0.65$ and rewrite the sum \eqref{TP4_a} in terms of a single index $m$ to mimic the form \eqref{a_decomp}, with the sequence $\{\nu_m\}_{m=1}^\infty$ ordered descendingly.

\begin{table}[t!]
\centering
\caption{Reference energies $||u^{\textrm{ref}}||_B$ for the four test problems TP.1--TP.4 presented in Section \ref{Sec:num_results}.}\label{Tab:reference_energies}
\setlength{\tabcolsep}{3pt}
\begin{tabu}{|[1.3pt] c |c|[1.3pt]}
\tabucline[1.3pt]{-}
Test Problem & Reference Energy $||u^{\textrm{ref}}||_B$ \\
\hline
  TP.1 & 1.50342524$\times10^{-1}$ \\
  TP.2 & 1.90117000$\times10^{-1}$ \\
  TP.3 & 1.94142000$\times10^{-1}$ \\
   TP.4  & 1.34570405$\times10^{-1}$ \\
\tabucline[1.3pt]{-}
\end{tabu}
\end{table}

\subsection{Experimental Setup}\label{Sec:exp_setup}
To begin, we select an appropriate set of finite element spaces $\mathbf{H}_1$. Since $D$ is square in all cases we choose a sequence $\bm{\mathcal{T}}$ of uniform meshes of square elements, with $\mathcal{T}_i$ representing a $2^i\times 2^i$ grid over $D$ (thus $\mathcal{T}_{i+1}$ represents a uniform refinement of $\mathcal{T}_i$) with element width $h(i)=2^{1-i}$ for TP.1 and $h(i)=2^{-i}$ for TP.2--TP.4. We then choose $\mathbf{H}_1$ to be the set of $\mathbb{Q}_1$ finite element spaces associated with $\bm{\mathcal{T}}$.  We initialise Algorithm \ref{Alg:1} with
\begin{align*}
J_P^0 = \{(0,0,\dots),(1,0,\dots)\},\qquad \bm\ell^0 = \{4,4\}\qquad (16\times 16 \textrm{ grids}).
\end{align*}

To compute the error estimator $\eta$ defined in Section \ref{Sec:a_post_error_para_diff}, the FEM spaces $\mathbf{H}_2 = \{H_2^\mu\}_{\mu\in J_P}$ are chosen to be broken $\mathbb{Q}_2$ spaces (see Figure \ref{Fig:H2detail}) defined with respect to the same meshes as the spaces $\mathbf{H}_1$, as described in Section \ref{Sec:spatial_apost}. Note that for this setup, if $a_0$ in \eqref{a_decomp} is a constant, we have $\gamma \leq \sqrt{5/11}$ in \eqref{B0_SCS}; c.f. Remark \ref{Rem:CBSrem2} and see \cite{CBS} for a proof. We also fix $\Delta_M = 5$ in the definition of $J_Q$ in \eqref{JQ_setup}.  Due to Galerkin orthogonality, the exact energy error $||u-u_X^k||_B$ at step $k$ admits the representation
\begin{align}\label{exact_error_decomp}
||u-u_X^k||_B = \left(||u||_B^2 - ||u_X^{k}||_B^2\right)^{\frac{1}{2}}.
\end{align}
To examine the \emph{effectivity index} $\theta^k = \eta^k/||u-u^k||_B$ we approximate $u$ in \eqref{exact_error_decomp} with an accurate `reference' solution $u^{\textrm{ref}}\in X^{\textrm{ref}}$. {The space $X^{\textrm{ref}}$ is generated by applying Algorithm \ref{Alg:1} with a much smaller error tolerance $\epsilon$ than the one used to generate $\eta^k,$ $k=1,\ldots, K$}. The reference energies $||u^{\textrm{ref}}||_B$  required for the approximation of \eqref{exact_error_decomp} are provided in Table \ref{Tab:reference_energies}. 

\begin{figure}[t!]
\centering
\includegraphics[width=0.95\textwidth]{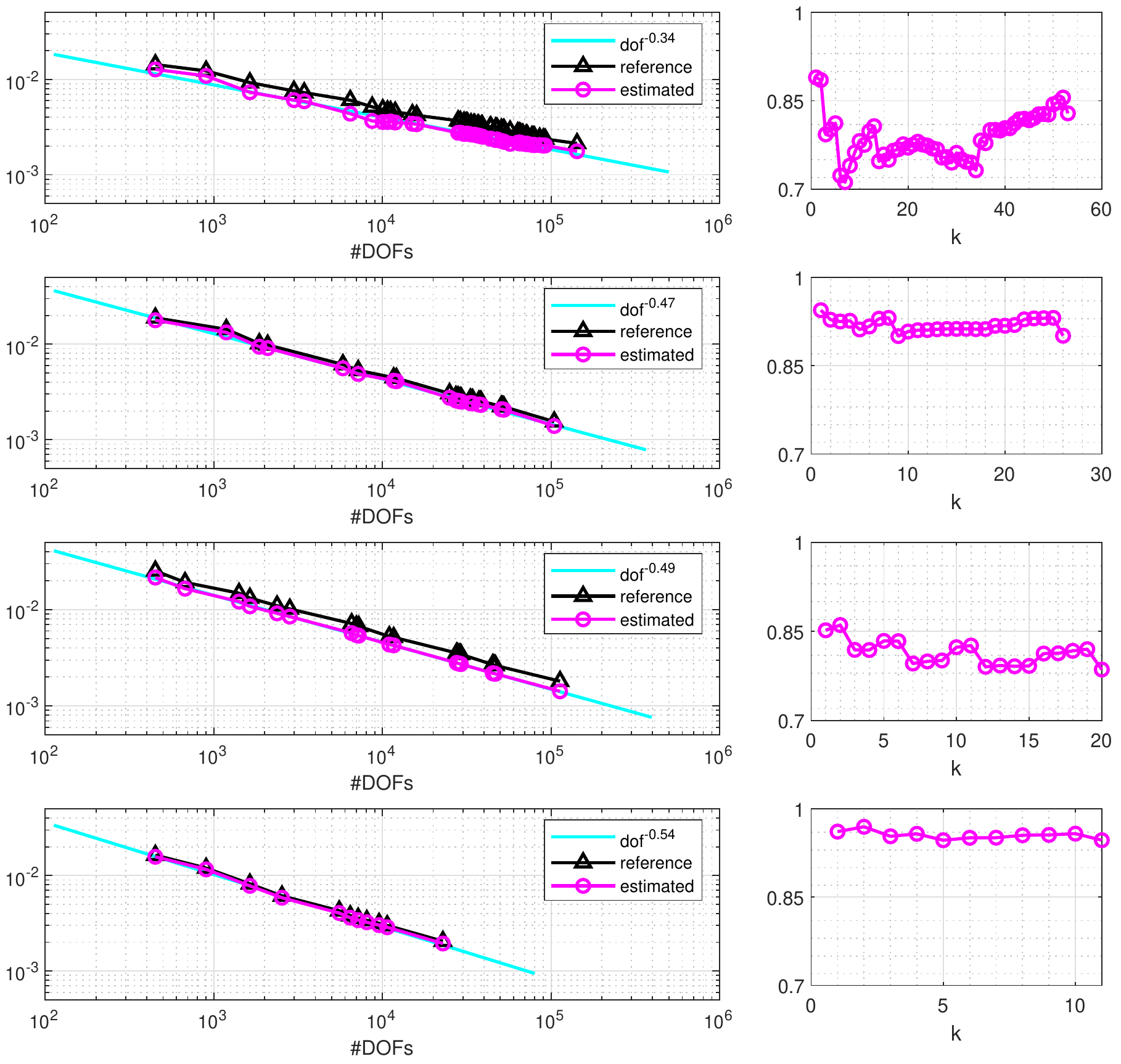}
\caption{Plots of the estimated errors $\eta^{k}$ versus number of degrees of freedom $N_{\textrm{dof}}$ (left) at steps $k=0,1,\ldots$ and effectivity indices $\theta^{k}$ (right) when solving TP.1--TP.4 (top-to-bottom) using Algorithms \ref{Alg:1} and \ref{Alg:2} (version 1). }
\label{Fig:exp1_1}
\end{figure}

\subsection{Experiment 1 (convergence rates)}\label{Sec:exp_1}
In our first experiment we solve test problems TP.1--TP.4 using Algorithms \ref{Alg:1} and \ref{Alg:2} (version 1) with tolerance $\epsilon = 2\times 10^{-3}$. In Figure \ref{Fig:exp1_1} we plot the evolution of the estimated error $\eta^k$ against $\textrm{dim}(X^k)$ (left plots) over each step of the iteration, as well as estimates of the effectivity indices $\theta^k$ (right plots). For test problems TP.2--TP.4, we observe that the estimated error behaves like $N_{\textrm{dof}}^{-0.5}$. Note that this is an improvement on the convergence rates obtained in \cite{MR3519560,MR3771406} for the same test problems, where single-level SGFEM spaces of the form \eqref{single_level_X} are employed. Due to our choice of FEM spaces $\mathbf{H}_1$ (bilinear approximation), and the spatial regularity of the solution, this is the optimal rate of convergence. That is, we achieve the rate afforded to the analogous parameter-free problem when employing $\mathbb{Q}_1$ approximation over uniform square meshes, and performing uniform mesh refinements. As proven in \cite{MR2728424,MR2763359,MR3091365}, the optimal achievable rate is a consequence of the fact that the sequence $\{||a_m||_\infty\}_{m=1}^\infty$ decays sufficiently quickly, and the error attributed to the choice of spatial discretisation dominates. Conversely, for test problem TP.1 the associated sequence $\{||a_m||_{\infty}\}_{m=1}^\infty$ decays too slowly, and the error attributed to the parametric part of the approximation dominates. For this reason, test problem TP.1 is particularly challenging. Nevertheless, for moderate error tolerances, our adaptive algorithm can tackle it efficiently. For all test problems considered, the effectivity indices are close to one, meaning that the error estimate is highly accurate.

\begin{table}[t!]
\centering
\caption{Number of solution modes assigned the same element width $h(\ell^\mu_K)$ (correspondong to a mesh level number $\ell_K^\mu$ in $\bm\ell^K$) for test problems TP.1--TP.4.}\label{Tab:exp1_1}
\setlength{\tabcolsep}{3pt}
\begin{tabu}{|[1.3pt] c |[1.3pt]c c c c c c|c | c|[1.3pt]}
\tabucline[1.3pt]{-}
Test Problem & $2^{-3}$ & $2^{-4}$ & $2^{-5}$ & $2^{-6}$ & $2^{-7}$ & $2^{-8}$ & $\textrm{card}(J_P^{K})$ & $M$ \\
\hline
TP.1 & 118 & 49 & 1 & 0 & 1 & 0 & 169 & 93\\
TP.2 & --  & 25 & 6 & 3 & 1 & 1 & 36  & 13\\
TP.3 & --  & 5  & 7 & 2 & 2 & 1 & 17  & 3 \\
TP.4 & --  & 17 & 3 & 0 & 1 & 0 & 21  & 8 \\
\tabucline[1.3pt]{-}
\end{tabu}
\end{table}

\begin{table}[t!]
\centering
\caption{A subset of $12$ multi-indices from the set $J_P^K$ generated by Algorithm \ref{Alg:1} and the associated element widths $h(\ell^\mu_K)$ assigned to those multi-indices at the final step for test problems TP.1--TP4.}\label{Tab:exp1_2}
\setlength{\tabcolsep}{3pt}
\begin{tabu}{|[1.3pt] c c|| c c|| c c|| c c|[1.3pt]}
\tabucline[1.3pt]{-}
 \multicolumn{2}{|[1.3pt]c||}{TP.1} & \multicolumn{2}{c||}{{TP.2}} & \multicolumn{2}{c||}{{TP.3}} & \multicolumn{2}{c|[1.3pt]}{{TP.4}}\\
\hline
 $\mu$ & $h(\ell^\mu_K)$ & $\mu$ & $h(\ell^\mu_K)$ & $\mu$ & $h(\ell^\mu_K)$ & $\mu$ & $h(\ell^\mu_K)$\\
\hline
 $(0\ 0\ 0\ 0\ 0\ 0\ 0\ 0\ 0\ 0)$			& $2^{-7}$ & $(0\ 0\ 0\ 0\ 0\ 0)$ 				  & $2^{-8}$ & $(0\ 0\ 0)$ 			& $2^{-8}$ & $(0\ 0\ 0\ 0\ 0\ 0)$ & $2^{-7}$\\
 $(\textbf{1}\ 0\ 0\ 0\ 0\ 0\ 0\ 0\ 0\ 0)$ & $2^{-5}$ & $(\textbf{1}\ 0\ 0\ 0\ 0\ 0)$ 		  & $2^{-7}$ & $(\textbf{1}\ 0\ 0)$ & $2^{-7}$ & $(\textbf{1}\ 0\ 0\ 0\ 0\ 0)$ & $2^{-5}$\\
 $(0\ 0\ \textbf{1}\ 0\ 0\ 0\ 0\ 0\ 0\ 0)$ & $2^{-4}$ & $(0\ 0\ \textbf{1}\ 0\ 0\ 0)$ 		  & $2^{-6}$ & $(\textbf{2}\ 0\ 0)$ & $2^{-7}$ & $(0\ 0\ \textbf{1}\ 0\ 0\ 0)$ & $2^{-5}$\\
 $(0\ \textbf{1}\ 0\ 0\ 0\ 0\ 0\ 0\ 0\ 0)$ & $2^{-4}$ & $(0\ \textbf{1}\ 0\ 0\ 0\ 0)$ 		  & $2^{-6}$ & $(\textbf{3}\ 0\ 0)$ & $2^{-6}$ & $(0\ \textbf{1}\ 0\ 0\ 0\ 0)$ & $2^{-5}$\\
 $(0\ 0\ 0\ 0\ 0\ \textbf{1}\ 0\ 0\ 0\ 0)$ & $2^{-4}$ & $(\textbf{2}\ 0\ 0\ 0\ 0\ 0)$ 		  & $2^{-6}$ & $(0\ \textbf{1}\ 0)$ & $2^{-5}$ & $(0\ 0\ 0\ \textbf{1}\ 0\ 0)$ & $2^{-4}$\\
 $(0\ 0\ 0\ 0\ \textbf{1}\ 0\ 0\ 0\ 0\ 0)$ & $2^{-4}$ & $(\textbf{1}\ \textbf{1}\ 0\ 0\ 0\ 0)$ & $2^{-5}$ & $(\textbf{4}\ 0\ 0)$	& $2^{-6}$ & $(\textbf{1}\ 0\ \textbf{1}\ 0\ 0\ 0)$ & $2^{-4}$\\
 $(0\ 0\ 0\ \textbf{1}\ 0\ 0\ 0\ 0\ 0\ 0)$ & $2^{-4}$ & $(0\ 0\ 0\ 0\ 0\ \textbf{1})$ 		  & $2^{-5}$ & $(\textbf{1}\ \textbf{1}\ 0)$ & $2^{-5}$ & $(\textbf{1}\ \textbf{1}\ 0\ 0\ 0\ 0)$ & $2^{-4}$\\
 $(\textbf{2}\ 0\ 0\ 0\ 0\ 0\ 0\ 0\ 0\ 0)$ & $2^{-3}$ & $(0\ 0\ 0\ 0\ \textbf{1}\ 0)$ 		  & $2^{-5}$ & $(\textbf{5}\ 0\ 0)$ & $2^{-5}$ & $(\textbf{2}\ 0\ 0\ 0\ 0\ 0)$ & $2^{-4}$\\
 $(0\ 0\ 0\ 0\ 0\ 0\ 0\ \textbf{1}\ 0\ 0)$ & $2^{-4}$ & $(0\ 0\ 0\ \textbf{1}\ 0\ 0)$ 		  & $2^{-5}$ & $(\textbf{2}\ \textbf{1}\ 0)$ & $2^{-5}$ & $(0\ 0\ 0\ 0\ 0\ \textbf{1})$ & $2^{-4}$\\
 $(0\ 0\ 0\ 0\ 0\ 0\ \textbf{1}\ 0\ 0\ 0)$ & $2^{-4}$ & $(\textbf{1}\ 0\ \textbf{1}\ 0\ 0\ 0)$ & $2^{-5}$ & $(0\ 0\ \textbf{1})$ & $2^{-5}$ & $(0\ 0\ 0\ 0\ \textbf{1}\ 0)$ & $2^{-4}$\\
 $(0\ 0\ 0\ 0\ 0\ 0\ 0\ 0\ 0\ \textbf{1})$ & $2^{-4}$ & $(\textbf{2}\ \textbf{1}\ 0\ 0\ 0\ 0)$ & $2^{-4}$ & $(\textbf{3}\ \textbf{1}\ 0)$ & $2^{-5}$ & $(\textbf{1}\ 0\ 0\ \textbf{1}\ 0\ 0)$ & $2^{-4}$\\
 $(0\ 0\ 0\ 0\ 0\ 0\ 0\ 0\ \textbf{1}\ 0)$ & $2^{-4}$ & $(\textbf{3}\ 0\ 0\ 0\ 0\ 0)$ 		  & $2^{-5}$ & $(\textbf{6}\ 0\ 0)$ & $2^{-5}$ & $(0\ \textbf{1}\ \textbf{1}\ 0\ 0\ 0)$ & $2^{-4}$\\
\tabucline[1.3pt]{-}
\end{tabu}
\end{table}

Figure \ref{Fig:exp1_1} provides no information about the structure of the multilevel SGFEM spaces $X^K$ constructed. To illustrate the qualitative differences between the four cases, in Table \ref{Tab:exp1_1} we record the number of activated parameters $M$, the cardinality of the final set $J_{P}^{K}$ and the number of multi-indices within that set that are assigned the same finite element space (i.e., the same mesh level number from the set $\bm\ell^K$). In each case, we observe that fine meshes are required to estimate very few solution modes (polynomial coefficients), whereas higher numbers of modes are assigned coarse meshes. This is reminiscent of multilevel sampling methods. While multilevel Monte Carlo and multilevel and multi-index stochastic collocation methods \cite{MR2835612,MR3033013,MR3416136,MR3502561,MR3579717} also typically require few deterministic PDE solves using fine finite element meshes and larger numbers using coarser meshes, there are some differences. Multilevel sampling methods typically require the number of parameters to be fixed a priori. We stress that our algorithm requires no sampling and learns which are the important parameters to activate as part of the solution process itself.  The decision about which meshes to use is based on an a rigorous a posteriori error estimate. For TP.1, we observe that many more parameters are activated ($M=93$) and the number of polynomials required ($\textrm{card}(J_P^{K})=169$) is much higher than in test problems TP.2--TP.4. This is due to the slow decay of the eigenvalues $\lambda_m$ in \eqref{TP1_a}. Although many more polynomials are needed in TP.1, the majority of the corresponding meshes are coarse. Conversely, test problem TP.3 has the lowest number of activated parameters ($M=3$) and requires the smallest number of polynomials ($\textrm{card}(J_P^{K})=17$). Compared to TP.1, however, a larger proportion of the meshes associated with the selected multi-indices are finer.  For TP.2, the number of activated parameters is higher than in TP.3, as expected.

In Table \ref{Tab:exp1_2} we display twelve of the multi-indices in the set $J_P^K$ that are selected by Algorithm \ref{Alg:1} for each test problem, as well as the associated element widths $h(\ell^\mu_K)$ assigned to those multi-indices, at the final step.  Note that it is not possible to list all the multi-indices generated for all four test problems. The twelve shown in each case are selected in the first few iterations.  For TP.1, these mostly correspond to univariate polynomials of degree one. In the early stages, Algorithm \ref{Alg:1} selects multi-indices that activate more terms in the expansion \eqref{TP1_a}, rather than multi-indices that correspond to polynomials of higher degree in the currently active parameters. Again, this is due to the slow decay of the $\lambda_m$ in \eqref{TP1_a}. In contrast, when solving  TP.3,  Algorithm \ref{Alg:1} first selects multi-indices that correspond to polynomials of higher degree in the currently active parameters, before activating new parameters.  For all test problems, the multi-indices that are selected in the early stages (corresponding to the most important solution modes, with respect to the energy error), are assigned the finest meshes. In particular, the \emph{mean} solution mode is the coefficient of the polynomial associated with $\mu=(0,0,\dots)$. This is always allocated the finest mesh.

\begin{table}[t!]
\centering
\caption{Solution times $T$ (in seconds) and adaptive step counts $K$ required to solve test problems TP.1--TP.4 using Algorithms \ref{Alg:1} and \ref{Alg:2} (versions 1 and 2) with various choices of the error tolerances $\epsilon$. The symbol `--' denotes that the estimated error at the previous step is already below the tolerance and the preceeding $T$ and $K$ are applicable.}\label{Tab:exp2_1}
\setlength{\tabcolsep}{3pt} 
\begin{tabu}{|[1.3pt]c|[1.3pt]cc|cc||cc|cc||cc|cc||cc|cc|[1.3pt]}
\tabucline[1.3pt]{-}
\multicolumn{1}{|[1.3pt]c|[1.3pt]}{\multirow{2}{*}{-}} & \multicolumn{4}{c||}{TP.1} & \multicolumn{4}{c||}{TP.2} & \multicolumn{4}{c||}{TP.3} & \multicolumn{4}{c|[1.3pt]}{TP.4}\\
\cline{2-17}
\multicolumn{1}{|[1.3pt]c|[1.3pt]}{} & \multicolumn{2}{c|}{ver. 1} & \multicolumn{2}{c||}{ver. 2} & \multicolumn{2}{c|}{ver. 1} & \multicolumn{2}{c||}{ver. 2} & \multicolumn{2}{c|}{ver. 1} & \multicolumn{2}{c||}{ver. 2}  & \multicolumn{2}{c|}{ver. 1} & \multicolumn{2}{c|[1.3pt]}{ver. 2}\\
\hline
$\epsilon$ & $T$ & $K$ & $T$ & $K$ & $T$ & $K$ & $T$ & $K$ & $T$ & $K$ & $T$ & $K$ & $T$ & $K$ & $T$ & $K$\\
\hline
$4.5\cdot10^{-3}$ & 2   & 6  & 2   & 6  & 1  & 7  & 5  & 6  & 1  & 10 & 1  & 7  & 1 & 5  & 2  & 5 \\
$3.0\cdot10^{-3}$ & 13  & 14 & 3   & 8  & 4  & 9  & -- & -- & 3  & 12 & 3  & 9  & 2 & 10 & -- & --\\
$1.5\cdot10^{-3}$ & 311 & 83 & 325 & 34 & 27 & 26 & 29 & 10 & 16 & 20 & 11 & 11 & 7 & 19 & 5  & 7 \\
\cline{2-5}
$9.0\cdot10^{-4}$ &\multicolumn{4}{c||}{} & 236 & 70 & 167 & 13 & 87 & 36 & 62 & 15 & 23 & 29 & 22 & 8\\
$7.5\cdot10^{-4}$ &\multicolumn{4}{c||}{\multirow{2}{*}{\underline{out of memory}}} & -- & -- & -- & -- & 100 & 38 & -- & -- & 36 & 38 & -- & --\\
$6.0\cdot10^{-4}$ &\multicolumn{4}{c||}{} & 881 & 147 &  --  & -- & 147 & 44 & 92 & 18 & 110 & 48 & 80  & 9\\
$4.5\cdot10^{-4}$ &\multicolumn{4}{c||}{} & 2197 & 177 & 1306 & 19 & 484 & 61 & 340 & 22 & 158 & 59 & 95 & 10\\
\tabucline[1.3pt]{-}
\end{tabu}
\end{table}

\subsection{Experiment 2 (timings)} We now investigate the computational efficiency of the new method. All computations were performed in MATLAB using new software developed from components of the S-IFISS toolbox \cite{SIFISS} on an Intel Core i7 4770k 3.50GHz CPU with 24GB of RAM. In Table \ref{Tab:exp2_1} we record timings ($T$) in seconds and the number of adaptive steps ($K$) taken by Algorithm \ref{Alg:1} (using both versions of Algorithm \ref{Alg:2} now), as we decrease the error tolerance $\epsilon$. We observe that for TP.2--TP.4, for smaller error tolerances, using version 2 of Algorithm \ref{Alg:2} results in a quicker solution time and a lower adaptive step count. The lower step count is due to the fact that the sets of multi-indices $\bar J_k$ that are produced by version 2 are usually richer than the ones produced by version 1. Note that because of this, a single step of version 2 is more expensive than a single step of version 1. Time savings are only made when enough steps are saved. We use the preconditioned conjugate gradient method with a mean--based preconditioner \cite{MR2491431} to solve \eqref{para_discrete_prob}. Fewer adaptive steps means that fewer SGFEM linear systems have to be solved and hence fewer matrix--vector products \eqref{multi_mat_vec} are required. For TP.1 with $\epsilon = 1.5\times10^{-3}$, the difference in step count between version 1 and 2 is not large enough for time savings to be made. We note also that asymptotically, both versions of Algorithm \ref{Alg:2} result in the same rates of convergence (illustrated by the blue lines in Figure \ref{Fig:exp1_1}). However, due to the larger associated sets $\bar J_k$, version 2 requires more adaptive steps before this rate is realised.

\begin{figure}[t!]
\centering
\includegraphics[width=0.93\textwidth]{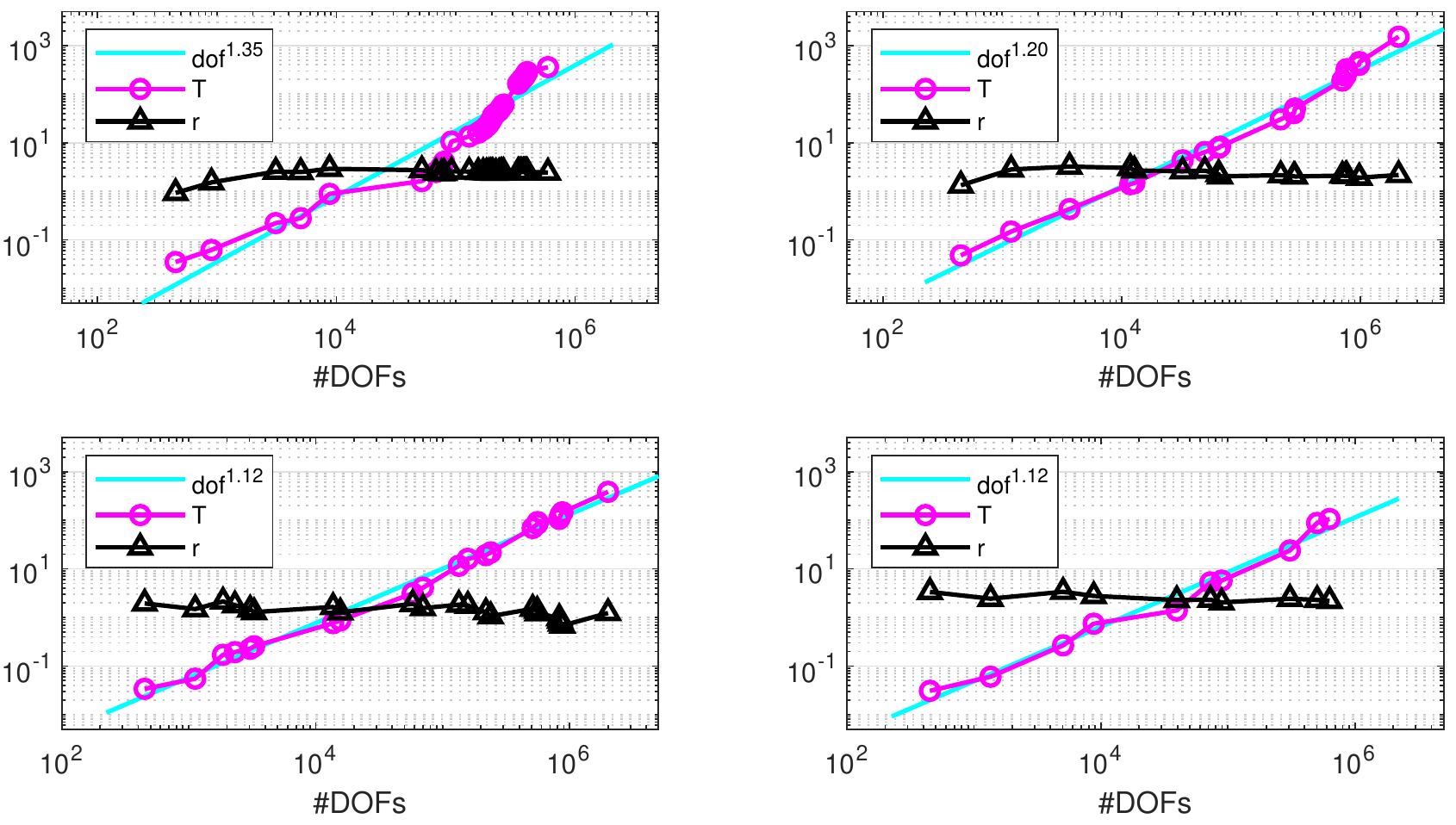}
\caption{Plots of the total computational time $T$ (round markers) in seconds accumulated over all refinement steps and the error estimation--solve time ratio $r$ (triangular markers), versus the number of degrees of freedom $N_{\textrm{dof}}$ when solving TP.1--TP.4 (left-to-right, top-to-bottom) using Algorithm \ref{Alg:1} with version 2 of Algorithm \ref{Alg:2}.}
\label{Fig:exp2_1}
\end{figure}

In Figure \ref{Fig:exp2_1} we plot the total computational time ($T$) against the the number of degrees of freedom ($N_{\textrm{dof}}$) when employing version 2 of Algorithm \ref{Alg:2}.   The total number of markers, each reflecting a single step of Algorithm \ref{Alg:1}, is equal to the value of $K$ corresponding to the smallest value of $\epsilon$ in Table \ref{Tab:exp2_1}. We observe that for all four test problems, the computational time behaves at most like $N_{\textrm{dof}}^{1.35}$. For TP.3 and TP.4, where $M$ is smaller, $T$ behaves almost linearly with respect to  $N_{\textrm{dof}}$. We also plot the ratio $r$ of the cumulative time taken to estimate the energy error (by executing the modules \texttt{COMPONENT\_SPATIAL\_ERRORS} and \texttt{PARAMETRIC\_SPATIAL\_ERRORS} in Algorithm \ref{Alg:1}) to the time taken to compute the SGFEM approximation $u_{X}$ (by executing the \texttt{SOLVE} module in Algorithm \ref{Alg:1}). We observe that $r$ does not grow with $N_\textrm{dof}$ (indeed, $1.2 < r < 2.5$ at the final step for all four problems). Hence, the cost of estimating the error is proportional to the cost of computing the SGFEM approximation itself. 

\section{Summary}
We presented a novel adaptive multilevel SGFEM algorithm for the numerical solution of elliptic PDEs with coefficients that depend on countably many parameters $y_{m}$ in an affine way. A key feature is the use of an implicit a posteriori error estimation strategy to drive the adaptive enrichment of the approximation space. We demonstrated how to extend the error estimation strategy used in \cite{MR3177362,MR3519560} to the new multilevel setting and described new ways to utilise the distinct components of the error estimator to determine how to best enrich the spaces associated with the spatial and parameter domains. Through numerical experiments we demonstrated that the error estimate is accurate and that the resulting adaptive algorithm achieves the optimal rate of convergence with respect to the dimension of the approximation space. That is, we achieve the convergence rate associated with the chosen finite element method for the associated parameter-free problems. Unlike other methods, our numerical scheme uses no marking or tuning parameters. Finally, we demonstrated that our multilevel algorithm is computationally efficient. Indeed, for some test problems (where the number $M$ of parameters that need to be activated is not too high), the solution time scales almost linearly with respect to the dimension of the approximation space.

\bibliography{bibdatabase}

\end{document}